\documentclass[preprint,12pt]{elsarticle}

\usepackage[T1]{fontenc}
\usepackage{amsmath,amssymb,latexsym,amsthm}
\usepackage{colortbl}

\usepackage{appendix}
\usepackage{url}
\usepackage{natbib}
\usepackage{graphicx}
\usepackage{calc}
\usepackage{xcolor}
\usepackage{subcaption}
\usepackage{diagbox}
\usepackage{booktabs}
\usepackage{textcomp}

\usepackage{wasysym}

\usepackage[normalem]{ulem}
\usepackage{soul}
\usepackage{array}
\usepackage{graphicx}
\usepackage{multirow}
\usepackage{textcomp}

\newcommand{\R}{\mathbb R}
\newcommand{\E}{\mathbb E}
\newcommand{\diag}{\mathrm{diag}}
\renewcommand{\hat}{\widehat}
\newcommand{\tr}{\mathrm{trace}}
\newcommand{\argmax}{\displaystyle \mathop{argmax}}

\graphicspath{{.}{./Figures/}}

\newbox\hautbox \setbox\hautbox=\hbox{\vphantom{\rule[-1cm]{0cm}{1.5cm}}}

\usepackage{amssymb}

\usepackage{soul}
\setstcolor{red}
\journal{Computational Statistics \& Data Analysis}

\begin{document}

\begin{frontmatter}

\title{Outlier detection in multivariate functional data through a contaminated mixture model}

\author[inst1,inst2]{Martial AMOVIN-ASSAGBA\corref{cor1}}
\cortext[cor1]{Corresponding author.\\
E-mail addresses: martial.amovin@masterk.com (M. AMOVIN-ASSAGBA), irene.gannaz@insa-lyon.fr (I. GANNAZ), julien.jacques@univ-lyon2.fr (J. JACQUES) }

\affiliation[inst1]{organization={Univ Lyon, Univ Lyon 2},
            addressline={ERIC UR3083}, 
            city={Lyon},
            country={ France}}
            
\affiliation[inst2]{organization={Arpege Master K},
            city={Saint-Priest},
            postcode={69800}, 
            country={France}}

\author[inst3]{Ir\`ene GANNAZ}


\affiliation[inst3]{organization={Univ Lyon, INSA Lyon, UJM, UCBL, ECL, ICJ},
            addressline={UMR5208}, 
            city={Villeurbanne},
            postcode={69621}, 
            country={France}}
            
\author[inst1]{Julien JACQUES}

\begin{abstract}

In an industrial context, the activity of sensors is recorded at a high frequency. A challenge is to automatically detect abnormal measurement behavior. Considering the sensor measures as functional data, the problem can be formulated as the detection of outliers in a multivariate functional data set. Due to the heterogeneity of this data set, the proposed contaminated mixture model both clusters the multivariate functional data into homogeneous groups and detects outliers. The main advantage of this procedure over its competitors is that it does not require to specify the proportion of outliers. Model inference is performed through an Expectation-Conditional Maximization algorithm, and the BIC is used to select the number of clusters. 
Numerical experiments on simulated data demonstrate the high performance achieved by the inference algorithm. In particular, the proposed model outperforms the competitors. Its application on the real data which motivated this study allows to correctly detect abnormal behaviors.

\end{abstract}


\end{frontmatter}

\section{Introduction}
Recent innovations in measuring devices and acquisition methods, combined with the intensive use of computer resources, enable the collection of discrete massive data on increasingly fine grids. 
The huge dimensions of the resulting data, leading to the well known curse of dimensionality issue, has motivated the development of functional data analysis methods. 
Thus, rather than considering an observation as high-dimensional vector, it is considered as a discrete observation of a functional random variable $X = X(t)_{\{t\in [a,b]\}}$, taking its values in a space of infinite dimension.
Working with these functional data makes it possible to take into account the intrinsic regularity of these sequential data.
Many techniques or methods for analyzing functional data have been developed along years. Refer to \citep{ferraty2003curves,ramsay2005springer} for a synthesis.

In the industrial application which motivated this work, the functional variables correspond to the measure of sensors implanted into a metrology device. Several sensors are implanted into the same device. Thus, one measure corresponds to $p$ curves, one per sensor. So we consider a multivariate functional random variable $X(t) = (X^{(1)}(t), \ldots, X^{(p)}(t)) \in \R^p$,$\; p \geq 2$.  The objective of this work is then to build a procedure allowing us to detect outliers in these multivariate functional data.  

\begin{figure}[!ht]
    \centering
    \includegraphics[scale=0.45]{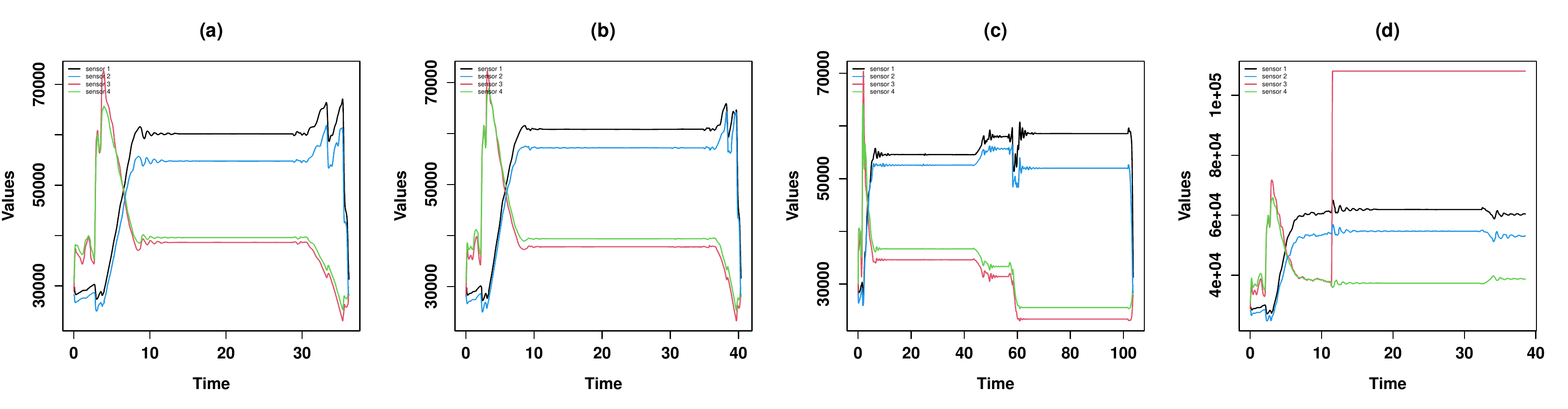}
    \caption{Example of normal sensor recordings (a, b) and of recordings presenting dysfunction (c, d).}
    \label{graph_data}
\end{figure}

The presence of outliers can significantly affect the results and conclusions of statistical analyses. In the industrial context, detecting outliers can be used to assess the quality of measurements from sensors or connected objects, and to identify abnormal or even faulty behavior. More precisely, it makes it possible to evaluate the reliability of measurements. 
Figure \ref{graph_data} presents a sample of normal curves and outliers from the industrial application which motivated this work, detailed in Section \ref{sec:real_data}. Data are measurements coming simultaneously from four sensors, each of the four curves on the same plot representing one sensor. The sensors work in pairs. They are used simultaneously to measure some quantity over a given period. A recording of the 4 sensors during the time period represents a datum.

Detecting outliers is challenging and many methods have been developed in the traditional (non-functional) multivariate case (see \emph{e.g.} \cite{agyemang2006comprehensive,chandola2009anomaly,hodge2004survey}). Some techniques are also being developed in the field of deep learning (see \emph{e.g.} \cite{chalapathy2019deep, braei2020anomaly}). In the functional data case, several methods have been recently introduced and are described below.
A first class of methods is based on the notion of depth. A depth function is a function $D(x,P)$ which measures the centrality or the closeness of a curve $x$ with respect to a probability distribution on the curves $P$. It gives a score which reflects the position of a curve relatively to the center of the data set. The higher the score, the more central the curve. A cutoff is set to specify whether the data are normal. Several depth functions have been defined by researchers to detect outliers in a functional setting.
In 2001, Fraiman and Muniz \citep{fraiman2001trimmed} introduced a depth function for univariate functional data and defined $\alpha-$trimmed means based on this notion ($\alpha$ the proportion of abnormal curves).  
Cuevas et al. \citep{cuevas2007robust}  defined a depth function based on the concept of mode and the random projection method. Febrero et al. \citep{febrero2008outlier} proposed a procedure for anomaly detection with depth functions which removes outliers at each iteration to avoid masking effects. Masking appears when some outliers mask the presence of others. Sun and Genton \citep{sun2011functional} constructed a functional boxplot to detect outliers.
In the multivariate case, Hubert et al. \citep{hubert2015multivariate, hubert2017multivariate} defined a new statistical distance named {\it bagdistance}, a measure derived from a depth function to rank the data. 
Dai and Genton \citep{dai2019directional} introduced the framework of a directional outlyingness for multivariate functional data. It is an effective measure of the variation in shape of the curves, which combines the direction with a statistical conventional depth.
A second class of methods is based on a distance measure named Dynamic Time Warping. It is used to find the similarity between the sequences of discrete observations of the functional data. The flexibility of the method enables the alignment of two time series (locally out of phase) in a non-linear manner \citep{sakoe1978dynamic}. By computing the similarity between time series two by two, it is possible to obtain a clustering by grouping the curves which are the closest, using any standard clustering algorithm \citep{giorgino2009computing, sarda2017comparing}. Since normal data look the same, they can form groups of normal data clusters. Data that are very far from normal data clusters can therefore constitute a cluster of atypical data.  It is worth noting that this method is particularly expensive in terms of execution. 
All the methods of these two first classes are based on a notion of distance, and require a threshold being set for this distance beyond which the observations are considered outliers.
A third class of methods consists of those interested in the notion of robust estimation. 
In 2005, \cite{garcia2005proposal} extended the trimmed kmeans method \citep{cuesta1997trimmed} to functional data. It is a robust method that allows data clustering, while excluding an \emph{a priori} known proportion of outliers.
Finally, a new method for anomaly detection in functional data has recently appeared, the Functional Isolation Forest \citep{staerman2019functional}, which does not use any notion of distance. It is based on an isolation technique  \citep{liu2008isolation} that separates abnormal data from others. Their method computes a path length which represents the number of cuts necessary to isolate a datum. The authors extend their model to the multivariate functional framework \citep{staerman2019functional}. The method requires an estimation of the number of outliers to define a threshold which will allow the specification of whether a datum is normal or not.

\begin{figure}[!ht]
    \centering
    \includegraphics[scale=0.35]{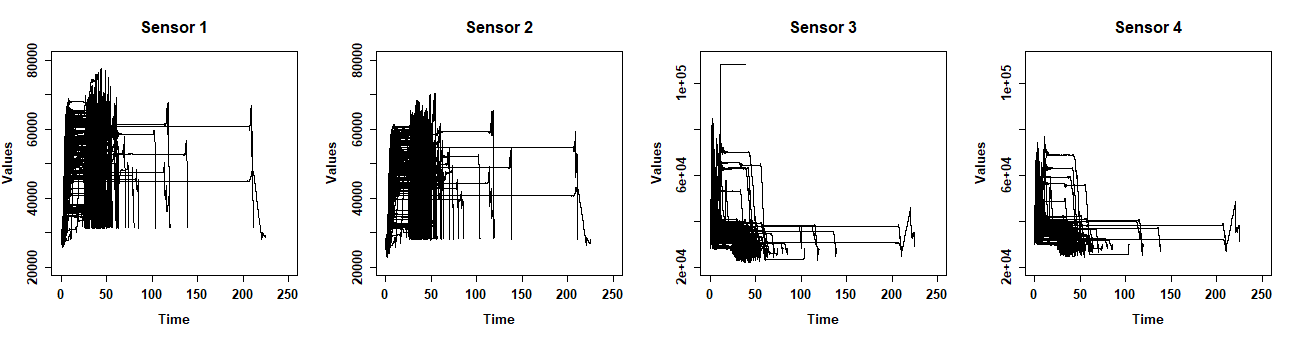}
    \caption{Representation of $148$ measurements, one plot per sensor.}
    \label{graph_data2}
\end{figure}

In the industrial application which motivated this work, the particularity of the measures is that they do not form an homogeneous data set, as illustrated by sensor on Figure \ref{graph_data2} which represents a sample of 148 measurements (one plot per sensor).
Indeed, this figure illustrates the differences both in the amplitude and in the duration of the measurements. These differences are mainly due to the differences between the entities which are measured.
Consequently, several clusters of measures co-exist, even in a normal operation without anomalies. 
The differences between the clusters can consequently mask the presence of outliers, or at least make their detection harder. 
Clustering functional data is a difficult task, since the data come from a space of infinite dimension. Several methods are proposed for overcoming this obstacle. A classification of them can be found in \cite{jacques2014functional}. The first approach is the two-stage method, which consists of reducing the dimension by approximating the data to elements from a finite dimensional space, and then using clustering algorithms from finite dimensional data \citep{abraham2003unsupervised,peng2008distance}.
Another method is the non-parametric approach, which uses usual clustering methods based on distances or dissimilarities specific to functional data \citep{ieva2011multivariate}. 
Finally, the model-based approach is one of a number of generative models that assume a probabilistic model either on the scores of a functional principal component analysis, or on the basis expansion coefficients of the approximations of the curves, in a basis of functions of finite dimension.
If there are many generative models in the univariate framework \citep{james2003clustering, heard2006quantitative,bouveyron2011model,jacques2014model}, only few are proposed in the multivariate framework. 
In 2013, \cite{jacques2013funclust} introduced a model-based clustering for multivariate functional data based on the Multivariate Functional Principal Component Analysis (MFPCA). Their model is based on cluster-specific Gaussian modelling of a given proportion of first FPCA scores.
In 2020, \cite{schmutz2020clustering} extended \cite{jacques2013funclust} by modeling all principal components whose estimated variance are non-null, which has the advantage of facilitating the model inference.  

In order to respond to the industrial motivation, the objective of the present paper is to be able to both cluster the observations into homogeneous groups and detect the presence of outliers. In the non-functional case, \cite{punzo2016parsimonious} introduced a contaminated multivariate Gaussian mixture model, which addresses our problem. According to their model, each component of the mixture is itself a mixture of normal and abnormal observations. Both have Gaussian distributions with a cluster-specific mean, but differ in their variance: the variance of the abnormal distribution is assumed to be the one of normal data multiplied by an inflation factor.
Such contaminated mixture models have also been considered in the non-Gaussian case \citep{tomarchio2020dichotomous, punzo2019new, punzo2018fitting}.


The present paper proposes an extension of \citep{punzo2016parsimonious} to the functional setting. For this, we will consider the model-based clustering procedure for multivariate functional data proposed by \citep{schmutz2020clustering}, and extend their mixture model to the contaminated case. In the resulting functional latent mixture model, normal and abnormal data are modeled into a functional subspace while using multivariate functional principal component analysis. The main advantage of the proposed model is that it does not require us to fix the proportion of outlier detection. A second advantage is that it allows us to consider model selection criteria for selecting the number of clusters.

The paper is organized as follows:  Section~\ref{sec:data} describes the data under study and introduces the notation and the specificity of functional data; Section~\ref{sec:proposed_model} presents the proposed contaminated mixture model for outlier detection in multivariate functional data; numerical experiments are detailed in Section~\ref{sec:num_exp}; and application to real data is discussed in Section~\ref{sec:real_data}.

\section{Multivariate functional data}
\label{sec:data}

We are seeking to define an automatic algorithm that simultaneously clusters multivariate functional data and detects outliers. 
Multivariate functional data are represented by a stochastic process $ X = \left\{X(t),\;t\in [a_1,a_2]\right\}, \; a_1, a_2 \in \mathbb R_{+}$. The realizations of the process $X$ belong to $L_2[a_1,a_2]$. Let  $ \{ X_1(t),X_2(t),\ldots,X_n(t) \}$ be an \emph{i.i.d.} sample of $ X$ where each copy $X_i$ of $X$ is a multivariate curve $X_i(t) = \left( X_i^1 (t),X_i^2(t), \ldots, X_i^p(t) \right)$ with $t\in [a_1,a_2]$. Let $x_i(t) = (x_i^1 (t),x_i^2(t),\ldots,x_i^p(t))$ denote the realizations of $X_i$, where $x_i^j(t)$ represents the measure of sensor $j$ during the period of time $t\in [a_1,a_2]$. In Section~\ref{sec:real_data} we will explain how we proceed in practice to reduce to such an interval for industrial application.

\subsection{Functional reconstruction}
In practice, the functions $x_i^j(t)$ are not totally observed, and we only have access to their values at certain discrete time points $a_1\leq t_1\leq t_2\leq\ldots\leq a_2$. There is a need to reconstruct the functional form of $x_i^j(t)$. One common approach is to assume that the curves $x_i^j(t)$ can be decomposed into a basis of functions $(\phi_b)_{1\leq b \leq B_j}$  in $L_2$, such that:
$$ x_i^j(t) = \sum_{b=1}^{B_j} c_{ib}^j \phi_b^j(t),\; t\in[a_1,a_2], \; 1\leq i \leq n,\; 1 \leq j \leq p,$$
where $(\phi_b^j)_{1\leq b \leq B_j}$ represents the basis functions of the $j^\text{th}$ component of the multivariate curve, $B_j$ the number of basis functions and $c_{ib}^j$ the coefficients.

Let the coefficients $(c_{ib}^j)$ be concatenated within a matrix $c$:  
$$ c =  (c_{B_1}^1, c_{B_2}^2, \ldots, c_{B_p}^p) \;\; \text{with} \;\; c_{B_j}^j = \begin{pmatrix}
c_{11}^j & \ldots & c_{1B_j}^j\\
c_{21}^j &  \ldots & c_{2B_j}^j\\
\vdots & \ddots &\vdots \\
c_{n1}^j & \ldots & c_{nB_j}^j \\
\end{pmatrix} . $$
$c_{B_j}^j$ represents the coefficients of the $j^\text{th}$ component of the multivariate curves.

Let the basis functions $(\phi_b)$ also be concatenated within a matrix $\Phi(t)$:
$$ \Phi(t) = \begin{pmatrix}
\Phi_{B_1}^1 (t) & 0_{B_2}  & \ldots & 0_{B_P}\\
0_{B_1} & \Phi_{B_2}^2(t)   &  \ldots & 0_{B_P}\\
\vdots & \vdots & \ddots &\vdots \\
0_{B_1} & 0_{B_2}  &  \ldots &   \Phi_{B_p}^p(t)   \\
\end{pmatrix}, 
$$
with $\Phi_{B_j}^j (t)  =  (\phi_{1}^j (t) , \phi_{2}^j (t), \ldots, \phi_{B_j}^j (t))$ and $0_\ell$ the null vector of size $\ell$.

With these notations, the set of multivariate curves $x(t)=(x_1(t),\ldots,x_n(t))$ can be written in the form:
\begin{equation}\label{eq:reconstruction}
    x(t) = c \Phi^\prime(t)\,,
\end{equation}
where the prime denotes the transpose operator.

The estimation of coefficients $c$ is generally performed by smoothing, and is amply described in \citep{ramsay2005springer}.
The choice of the basis as well as the number of basis functions strongly depends on the nature of data, and is done empirically in the unsupervised context. For example, Fourier basis are often used in the case of data with periodic pattern, spline basis for smooth curves and wavelets basis for very irregular curves. Depending on the data complexity, the number of basic functions can be increased or reduced.

\subsection{Probabilistic modeling}

The notion of probability density for functional variables is not well defined. In \cite{delaigle2010defining}, it is proved that it can be approximated with the probability density of the functional principal component scores (FPCA scores, \cite{ramsay2005springer}). Under the assumption \eqref{eq:reconstruction} of basis expansion decomposition, these FPCA scores are obtained directly from a PCA of the coefficient $c$ using a metric defined by the scalar product between the basis functions $\Phi$. Consequently, model-based approaches for functional data consider probabilistic distribution for either the FPCA scores \citep{jacques2013funclust, jacques2014model} or the basis expansion coefficients \citep{bouveyron2011model,schmutz2020clustering}, both being equivalent.

\section{A model for clustering contaminated functional data}
\label{sec:proposed_model}

\subsection{The Contaminated-funHDDC model}
\label{sec:model}

The aim of our work is to classify the observed multivariate curves $x_1, x_2, \ldots,x_n$ into $K$ clusters and to detect the presence of outliers. 
Let $z_i, 1\leq i \leq n$, where $z_i = (z_{i1},\dots,z_{iK})^\prime$, be a latent variable associated to each observation $x_i$, such that $z_{ik} = 1$ if the observation $x_i$ belongs to the cluster $k$ and $z_{ik} = 0$ otherwise. Let us similarly define $v_i, 1\leq i \leq n$, where $v_i = (v_{i1},\dots,v_{iK})^\prime$, such that $v_{ik} = 1$ if $x_i$ in group $k$ is a normal observation and $v_{ik} = 0$ if it is an outlier.
Hence, knowing that each curve $x_i$ can be identified by its basis expansion coefficient $c_i$, the completed-data can be expressed as $\mathcal{D} = \{c_i,z_i,v_i\}_{i=1}^n.$ 

Inspired by \cite{punzo2016parsimonious}, we propose the Contaminated-funHDDC (C-funHDDC) model:
\begin{equation}\label{eq:model}
    g(c_i,\theta) = \sum_{k=1}^K \pi_k \left[ \beta_k f(c_i,\mu_k,\Sigma_k) + (1-\beta_k)f(c_i,\mu_k,\eta_k\Sigma_k) \right],
\end{equation}
where $\theta = \{ \pi_k, \beta_k, \mu_k, \Sigma_k, \eta_k \}_{k=1}^K$ denotes all the parameters that need to be estimated. For all $k\in\{1,\dots, K\}$, $\pi_k \in (0,1]$, $\sum_{l=1}^K \pi_l =1$, $\eta_k >1$ is a covariance inflation factor, $\beta_k$ the proportion of normal data, $f(\cdot,\mu_k,\Sigma_k)$ the density of a Gaussian distribution of mean $\mu_k$ and covariance matrix $\Sigma_k$, $\beta_k \in [0,1]$ the proportion of normal observations and $g$  the density of the mixture of multivariate contaminated Gaussian distributions. 

The mixture-likelihood function is defined as 
\[
    L(c,\theta) =  \prod_{i=1}^n g(c_i, \theta) = \prod_{i=1}^n  \sum_{k=1}^K \pi_k \left[ \beta_k f(c_i,\mu_k,\Sigma_k) + (1-\beta_k)f(c_i,\mu_k,\eta_k\Sigma_k)\right].
\]

Depending on the regularity of the observed curves, the number of basis functions used for reconstruction can be large.
In order to introduce parsimony in the modeling, let us assume that the actual stochastic process associated with the $k^\text{th}$ cluster can be described in a low-dimensional functional latent subspace $\E_k[a_1,a_2]$ of $L_2[a_1,a_2]$ with dimension $d_k \leq  B= \sum_{j=1}^p B_j$, for every $k\in \{1,\ldots,K\}$. 
The curves of cluster $k$ can be expressed in a group-specific basis obtained from $\left\{\phi_b^j, \; j \in \{ 1,\ldots,p \},\; b\in \{1,\ldots,B\}\right\}$ by a linear transformation using a multivariate functional principal component analysis. That is, we introduce a basis $\left\{\varphi_{kb},\; b\in \{1,\ldots,B\}\right\}$ specific to cluster $k$ as \[\varphi_{kb} (t) = \sum_{l=1}^B q_{kbl}\phi_l(t), \;1 \leq b \leq B,\] with $q_{kbl}$ the coefficients of the eigenfunctions expressed in the initial base $\Phi$. 

Let $Q_k = (q_{kbl})_{\{1\leq b,l\leq B\}} $ an orthogonal matrix of size $B\times B$. We suppose that the $d_k$ first eigenfunctions contain main information of the MFPCA of cluster $k$. We therefore split $Q_k$ into two parts: $Q_k = [U_k,V_k]$ where $U_k$ is a matrix of size $B\times d_k$ and $V_k$ of size $B\times(B-d_k)$ and $U_k$ and $V_k$ are such that $U_k^\prime U_k = I_{d_k}$, $V_k^\prime V_k = I_{B-d_k}$ and $U_k^\prime V_k = 0$. This step will enable us to reduce the dimension of the modeling, since only first $d_k$ components will need a fine modeling.

The relationship between the original space of $c_{i}$ and the low-dimensional functional subspace is generated by $U_k$:
\begin{equation}\label{eq:coefficients_subspace}
    c_{i}  =   W^{-\frac{1}{2}} U_k \delta_{i}  + \epsilon_{i},
\end{equation}
where $W = \int_a^b \Phi^{\prime}(t) \Phi(t) \mathrm{d}t $ denotes the symmetric block-diagonal $ B \times B $ matrix of the inner products between the basis functions. $\epsilon_i$ is a noise following a Gaussian distribution conditionally to $z_{ik}=1$, 
\begin{equation*}
    \epsilon_{i}\mid z_{ik}=1 \sim \mathcal{N}(0,\Lambda_k),
\end{equation*}
where the covariance matrix $\Lambda_k$ will be defined just above.

Based on  \citep{delaigle2010defining}, it is possible to make distribution hypotheses on scores $\delta_{i}^\prime$. 
Let $n_k = \sum_{i=1}^n z_{ik}$ be the number of curves which belong to cluster $k$. 
Conditionally to $z_{ik}=1$, let assume that the principal component scores of the $n_k$ curves, $\delta_{i}^\prime, 1\leq i \leq n_k$, follow a Gaussian distribution. That is,
\begin{equation}
    \delta_{i} \mid z_{ik}=1, v_{ik}=1 \sim  \mathcal{N}(m_k, \Delta_k) \,,\;\; 
    \delta_{i} \mid z_{ik}=1, v_{ik}=0 \sim  \mathcal{N}(m_k, \sqrt{\eta_k}\Delta_k),
\end{equation}
with $m_k \in \R^{d_k}$, $\eta_k >1$ and $\Delta_k = \diag(a_{k1},a_{k2},\ldots,a_{k d_k}) $.

Hence the basis expansion coefficients satisfy 
\begin{equation}
   c_{i} \mid z_{ik}=1, v_{ik}=1 \sim \mathcal{N}\left(\mu_k , \Sigma_k  \right) \;,\;\;
   c_{i} \mid z_{ik}=1, v_{ik}=0 \sim \mathcal{N}\left(\mu_k , \eta_k\Sigma_k  \right),
\end{equation}
with $$ \mu_k =  U_k m_k, \;\;  \Sigma_k = W^{-\frac{1}{2}} U_k  \Delta_k U_k^\prime   W^{-\frac{1}{2}} + \Lambda_k\,.$$

Finally, let us assume that the noise covariance matrix $\Lambda_k$ is such that 
\begin{eqnarray}
\label{eq:HDDC}
 R_k= Q_k^\prime W^{\frac{1}{2}} \Sigma_k W^{\frac{1}{2}} Q_k  = \diag(a_{k1},\ldots,a_{kd_k},b_k,\ldots,b_k),
\end{eqnarray}
with $ a_{k1} >  a_{k2} > \ldots >  a_{kd_k} >  b_{k}$. In view of these notations, it is possible to conclude that  the variance of the actual data of the $k^\text{th}$ cluster is modeled by $ a_{k1}, a_{k2}, \ldots, a_{kd_k}$, while the remaining ones are considered as noise components and modeled by a unique parameter $b_{k}$. $R_k$ represents the variance of the data in the space generated by the eigenfunctions $\varphi_{kb}$.

\paragraph{Model identifiability}
We investigate in this paragraph the identifiability of C-funHDDC.
In \citep{yakowitz1968identifiability}, the identifiability of multivariate Gaussian mixture model is proved. 
However, according to \citep{di2007mixture}, in absence of any constraint, a mixture of mixtures is not identifiable. Punzo and McNicholas \citep{punzo2016contaminatedmixt} defined conditions allowing the identifiability of multivariate contaminated Gaussian mixture models which are actually a mixture of mixtures. As shown in \citep{punzo2016contaminatedmixt}, it is identifiable provided that, the $K$ Gaussian distributions representing the normal distributions in the mixture have distinct means and/or non-proportional covariance matrices. Let $k_1$ and $k_2$ denote the indexes of two different clusters:  $$ \text{if}\;\;  k_1 \neq k_2 \text{ then } \| \mu_{k_1} - \mu_{k_2} \|_2^2 + \| \Sigma_{k_1} - \kappa\Sigma_{k_2} \|_2^2 \neq 0 \;, \;\; \forall \kappa > 0, $$ 
where $\|\cdot\|_2$ is the Frobenius norm, then a finite mixture of contaminated Gaussian distributions is identifiable.
Since C-funHDDC is an extension of \citep{punzo2016contaminatedmixt} to a functional framework, then the identifiability of model \eqref{eq:model} in the space of the coefficients deduces directly from the identifiability of the mixtures of multivariate contaminated Gaussian distribution \citep{punzo2016contaminatedmixt}.

\paragraph{Model complexity}

Here the complexity of C-funHDDC, \emph{i.e.} the number of parameters to estimate is explored and compared with other models. For the numbers of parameters, let us define
$h = KB+K-1$, the number of parameters
required for the estimation of means and proportions, $w = \sum_{k=1}^K d_k[B-(d_k+1)/2]$ the number of parameters required for the estimation of orientation matrices $Q_k$, and $D = K + \sum_{k=1}^K d_k$, the number of parameters required for the estimation of $b_k$ and $a_k$. Due to the presence of outliers, we have $2K$ additional parameters which represent the number of parameters required for the estimation of $\{\beta_k\}_{k=1}^K$ and $\{\eta_k\}_{k=1}^K$. In total, the number of parameters to be estimated is: $h + w + D + 2K$.
The complexity of C-funHDDC is compared to the complexity of a usual multivariate Gaussian Mixture Model applied directly on the discretized data (discrete-GMM), and the C-funHDDC model without the parsimony assumption (\ref{eq:HDDC}) (functional-GMM). In the non-functional framework, the discretized data are considered as a high-dimensional vector, $x_i \in \R^{Tp}$,  where $T$ is the number of observation points $t_i$.
Let us see what happens with the real data. The observations per sensor are of the order of $T = 1000$ and we have $p=4$ sensors. The choices of $d_k=10$ and $K=3$ are explained in Section~\ref{sec:real_data}. The number of basis function ($B_j, 1\leq j \leq p$) per component is fixed at $25$. This implies $B = 100$. 
Table~\ref{table:model_complexity} presents the complexity of these models.
\begin{table}[!ht]
\centering
    \begin{tabular}{llr}
    \toprule
    Model& Parameter numbers & Param. nb. for  $T=1000$, \\
    &&$K=3$, $p=4$, $d_k=10$,$B=100$ \\
    \midrule
    C-funHDDC & $h + w + D + 2K$ & $ 3176 $  \\ 
    discrete-GMM & $(K-1)+T\{Kp +K[p(Tp+1)/2]\}$ &  $ 24018002 $  \\ 
    functional-GMM & $ (KB-1)+ 4K +K(B-1)(1+\frac{B}{2})  $ &  $ 15458 $ \\
    \bottomrule
    \end{tabular}
    \caption{Model complexity}
    \label{table:model_complexity}
\end{table}

The complexity of d-GMM is huge, illustrating that the functional approach is necessary for such data.  The functional approach should be favoured.

\subsection{Inference of the C-funHDDC model}
We have described the proposed latent block model. Now, we propose an algorithm for estimating the model parameters.

In model-based clustering, a classical way to estimate the parameters is the EM algorithm introduced by Dempster et al. \citep{dempster1977maximum}. It looks for an estimator maximizing the completed log-likelihood of the model by successive iterations of two steps E (Expectation) and M (Maximization).

The completed log-likelihood can be expressed as: 
$$\ell_c(\theta ) =  \ell_{c_1}(\theta) + \ell_{c_2} (\theta)  + \ell_{c_3} (\theta ),$$
with 
\begin{gather*}
\ell_{c_1}(\theta)  =  \sum_{i=1}^n \sum_{k=1}^K z_{ik} \log \pi_k,\quad 
\ell_{c_2} (\theta) = \sum_{i=1}^n \sum_{k=1}^K z_{ik} [ v_{ik} \log \beta_k + (1-v_{ik}) \log (1-\beta_k) ],\\
\text{and} \;\; \ell_{c_3} (\theta) = \sum_{i=1}^n \sum_{k=1}^K z_{ik} [ v_{ik} \log g(c_i,\mu_k,\Sigma_k) + (1-v_{ik}) \log g(c_i,\mu_k, \eta_k \Sigma_k) ].
\end{gather*}

To fit our model, we use a generalized EM algorithm called an Expectation-Conditional-Maximization (ECM) algorithm \citep{meng1993maximum}. It is an algorithm regularly used for maximum likelihood estimation when there is contamination or when data are missing \citep{browne2011model}. ECM algorithm differs from EM algorithm at the maximization step. Each M-step of EM algorithm, is replaced by two Conditional Maximization (CM) steps. They arise from the partition $\theta = \{\theta_1,\theta_2\}$ with $ \theta_1 = \{\pi_k, \beta_k, \mu_k,a_{kl}, b_k  \}_{k=1}^K$ and $\theta_2 = \{ \eta_k \}_{k=1}^K$. At the first CM step (CM1), we estimate the parameters $\theta_1$ conditionally to $\theta_2$, and at the second CM step (CM2), we estimate $\theta_2$ conditionally to $\theta_1$. As described below, the estimators of $\mu_k$ and $\Sigma_k$ depend on $\eta_k$ and vice versa. So it is necessary to estimate $\eta_k$ apart and fix this parameter before estimating $\mu_k$ and $\Sigma_k$ and vice versa.

\paragraph{Algorithm} The ECM algorithm starts by initializing the latent variables $z_{ik}^{(0)}$, $v_{ik}^{(0)}$ for each observation and the parameter $ \theta_2^{(0)} = \eta_k^{(0)}$. Then a first execution of CM1 step and CM2 step is undertaken in order to provide the information of the different groups necessary for computing the conditional expectation of E step (described below). Finally, we alternate E step, CM1 step and CM2 step until convergence ($\ell_c(\theta^{(m)})-\ell_c(\theta^{(m-1)})$ is less than a given threshold $\varepsilon$ or a maximum number of iterations is reached ($\theta^{(m)}$ is the value of $\theta$ at iteration $m$ of the ECM algorithm). We choose $\varepsilon = 0.0001$.

During the ECM algorithm $K$ and $d_k$ are considered as fixed hyper-parameters. The procedures for estimating these parameters are respectively described in subsection \ref{subsec:estimation_of_K} and \ref{subsec:estimation_of_dk}.

Now let us give the expression of parameters obtained at iteration $m$.

\textbf{E step (Expectation).} This step needs to compute $Q(\theta \mid \theta^{(m)}) = \mathbb{E} \left( \ell_c (\theta) \mid \theta^{(m-1)} \right)$. For that, we have to calculate $\E (z_{ik}\mid c_i,\theta^{(m)})$ and $\E (v_{ik}\mid c_i,\theta^{(m)})$ whose expressions are given below.
The posterior probability that the observation  $i$  belongs to the cluster $k$ of the mixture, can be expressed as:
\begin{align*}
    t_{ik}^{(m)} &:= \E (z_{ik}\mid c_i,\theta^{(m)})
    = \frac{\pi_k^{(m)} \left[ \beta_k^{(m)} f(c_i,\mu_k^{(m)},\Sigma_k^{(m)}) + (1-\beta_k^{(m)})f(c_i,\mu_k^{(m)},\eta_k^{(m)}\Sigma_k^{(m)}) \right]}{\sum_{l=1}^K \pi_l^{(m)} \left[ \beta_l^{(m)} f(c_i,\mu_l^{(m)},\Sigma_l^{(m)}) + (1-\beta_l^{(m)})f(c_i,\mu_l^{(m)},\eta_l^{(m)}\Sigma_l^{(m)})\right]}\,.
\end{align*}
The probability of an observation $i$ being normal or abnormal, given  that it belongs to cluster $k$, can be expressed as:
\begin{align*}
    s_{ik}^{(m)} := \E (v_{ik}\mid c_i,\theta^{(m)}) = \frac{\beta_k^{(m)} f(c_i,\mu_k^{(m)},\Sigma_k^{(m)})}{ \beta_k^{(m)} f(c_i,\mu_k^{(m)},\Sigma_k^{(m)}) + (1-\beta_k^{(m)})f(c_i,\mu_k^{(m)},\eta_k^{(m)}\Sigma_k^{(m)})}\,.
\end{align*}

\textbf{CM1 step (First conditional maximization).}
At this step, we fix $\theta_2$ at $\theta_2^{(m)}$ and we compute $\theta_1^{(m+1)}$ as the value of $\theta_1$  that maximizes $Q(\theta\mid\theta^{(m)})$. The estimation of the parameters is:
\begin{itemize}

\item $\pi_k^{(m+1)} = \frac{\sum_{i=1}^n t_{ik}^{(m)}}{n}$,

\item $ \beta_k^{(m+1)} = \frac{\sum_{i=1}^n t_{ik}^{(m)} s_{ik}^{(m)}}{\sum_{i=1}^n t_{ik}^{(m)}}$,

\item $\mu_k^{(m+1)} =  \frac{1}{\sum_{i=1}^n t_{ik}^{(m)}\left( s_{ik}^{(m)} + \frac{1-s_{ik}^{(m)}}{\eta_k^{(m)}} \right)} \sum_{i=1}^n t_{ik}^{(m)}\left( s_{ik}^{(m)} + \frac{1-s_{ik}^{(m)}}{\eta_k^{(m)}} \right) c_i$.

For the variance parameters, let us define the current sample covariance matrix of the cluster $k$ by:
\[
    H_k^{(m)} = \frac{1}{\gamma_k^{(m)}} \left( \sum_{i=1}^n   t_{ik}^{(m)}  \left( s_{ik}^{(m)} + \frac{1-s_{ik}^{(m)}}{\eta_k^{(m)}} \right) (c_i-\mu_k^{(m+1)})^\prime (c_i - \mu_k^{(m+1)}) \right),
\]
with $\gamma_k^{(m)} = \sum_{i=1}^n t_{ik}^{(m)}$. 

\item the $d_{kl}$ first columns of the matrix $Q_k$ are updated by the eigenfunctions
coefficients associated with the largest eigenvalues of $W^{\frac{1}{2}} H_k^{(m)} W^{\frac{1}{2}}$.

\item the variance parameters $a_{kl}^{(m+1)}$ are updated by the $d_k$ largest eigenvalues $ \lambda_{kl}^{(m)}$
of $W^{\frac{1}{2}} H_k^{(m)} W^{\frac{1}{2}}$  with $\lambda_{kl}^{(m)} = q_{kl}^\prime W^{\frac{1}{2}} H_k^{(m)} W^{\frac{1}{2}} q_{kl}$ and $q_{kl}$ the $l$th column of $Q_k$.  

\item the variance parameters $b_{k}$ are updated by $ \frac{1}{B-d_k}\left( \tr(W^{\frac{1}{2}} H_k^{(m)} W^{\frac{1}{2}}) - \sum_{l=1}^{d_k} \lambda_{kl}^{(m)}   \right)$.  

\end{itemize}

\textbf{CM2 step (Second conditional maximization).}
At this step, we fix $\theta_1$ at $\theta_1^{(m+1)}$ and we compute $\theta_2^{(m+1)}$ as the value of $\theta_2$  that maximizes $Q(\theta\mid\theta^{(m)})$. The expression of the estimation of the parameter $\eta_k$, is 
 $$\eta_k^{(m+1)} = \max\left\{ 1, \frac{\sum_{i=1}^n t_{ik}^{(m)}(1-s_{ik}^{(m)}) (c_i - \mu_k^{(m+1)}) ^\prime\Sigma_k^{-1(m+1)}(c_i - \mu_k^{(m+1)})}{B(\gamma_k^{(m)} - \gamma_{k_g}^{(m)})}  \right\},$$
 where  $\gamma_{k_g}^{(m)} = \sum_{i=1}^n t_{ik}^{(m)} s_{ik}^{(m)}$ is the current number of normal data in the cluster $k$.

The proof of the results above can be found in ~\ref{sec:sample:appendix}.

\paragraph{Classification} 
In order to define the group to which each observation belongs, we use the Maximum {\it a  posteriori} (MAP) classification. Let $i\in\{1,\dots,n\}$ and consider $\hat{z_i}$ and $\hat{v_i}$ the respective estimations of $z_i$ and $v_i$. 

For estimations $\{\hat{z}_{ik}, \;k=1,\dots,K\}$ of latent variables $\{z_{ik},\;k=1,\dots,K\}$, the MAP estimator is
$$
\text{MAP}(\hat{z}_{ik}) = \begin{cases}
        1 & \mbox{if } k=\argmax_{l=1,\dots,K} \hat{z}_{il},\\
        0 & \mbox{otherwise.}
    \end{cases}
$$

A datum $i$ belonging to cluster $k$ would be defined as normal if the expected value $\hat{v}_{ik}$ of $v_{ik}$ arising by the model is greater than $0.5$, and abnormal otherwise.

\subsection{Initialization}
The choice of the initialization method in EM or a generalized EM algorithm is a major issue in the speed of convergence and the performance of the clustering \citep{punzo2016parsimonious}. The parameters that are determined in the initialization are the latent variables $\{z_{ik}\}$ and $\{v_{ik}\}$ and the parameters $\{\eta_k\}$.

For $z_i,\; 1\leq i \leq n$, we propose three ways to initialize the algorithm. The first way is to randomly initialize the algorithm as in \citep{schmutz2020clustering, bouveyron2011model}. Another popular and common practice is the well known  \emph{kmeans} method \citep{hartigan1975clustering}.
Since the data are contaminated, these first two approaches do not give satisfactory results because they are not robust enough.
Another approach, which is the one we recommend to be used, is the {\it trimmed} approach. This method is a robust kmeans insensitive to outliers introduced by \citep{cuesta1997trimmed}. The main idea of {\it trimmed} is to calculate centroids by excluding a certain proportion of the data. We choose to take a large proportion of outliers at this step to ensure that we are not contaminated by their presence. This proportion is fixed at $0.2$. This will ensure our model does not depend on this parameter.

The ECM algorithm is a generalized EM algorithm \citep{meng1993maximum}, so it has the same log-likelihood growth properties. Hence it converges to a local maximum. Initialization of parameters $\{z_i,\; 1\leq i \leq n\}$ appears as the most critical. Therefore, to ensure a better convergence to a local maximum, we apply a "multiple initialization". That is, the algorithm is no longer executed once but a given number of times, denoted $nb\_init$, with different initialization values of $\{z_i,\; 1\leq i \leq n\}$. We considered $nb\_init = 10$ in our applications. Next, we retain the best of the $nb\_init$ results given by ECM algorithm using the BIC. 
We also initialize the latent variables $v_i,\; 1\leq i \leq n$ which specify if the corresponding observation is normal or not. For all $i \in \{1,\ldots, n\}$ and $k \in \{1,\ldots, K\}$ we set $v_{ik} = 0.99$ when $z_{ik}=1$ and $v_{ik} = 0$  when $z_{ik} = 0$. 

For each cluster the parameter $\eta_k, 1\leq k \leq K$ is initialized at $1$.

\subsection{Estimation of the intrinsic dimension $d_k$}
\label{subsec:estimation_of_dk}

Recall that for all $k=1,\dots, K$, the hyper-parameter $d_k$ (section \ref{sec:model}), is the intrinsic dimension of the cluster $k$. It needs to be estimated. We propose two methods of estimation.

The first method is the grid search method. $d_k$ is obtained by considering a grid of positive integers and choosing the best with the BIC.
The second way is to use the Cattell scree test \citep{cattell1966scree}. This test is often used in statistics to detect the \emph{bends} indicating a change of structure on the curve representing the eigenvalues.
As in \citep{schmutz2020clustering} we choose the value for which the subsequent eigenvalue differences are less than a threshold defined by the user. This value varies during iterations as in different groups. It depends on the MFPCA values obtained for each cluster.

\subsection{Choice of the number of clusters}
\label{subsec:estimation_of_K}
At this point, we are interested in selecting the number of cluster $K$. The choice of this hyper-parameter can be considered as a model selection problem. Many criteria are used in model selection. We can quote {\it e.g.} the ICL criterion \citep{biernacki2000assessing} or the Slope Heuristic criterion \citep{birge2007minimal} (see {\it e.g.} Table~2 of \citep{punzo2016contaminatedmixt}). In 1974, Akaike \citep{akaike1974new} proposed to maximize the likelihood function separately for each model, and then choose the model that maximizes the criterion $$ AIC =  l(\hat{\theta}) -\xi,$$  where  $\xi$ is the number of free parameters (Table~\ref{table:model_complexity}). $l(\hat{\theta})$ is the maximum log-likelihood value.

Another criterion was defined by Schwarz \citep{schwarz1978estimating} and expressed as follows:
$$ BIC = l(\hat{\theta}) - \frac{\xi}{2} \log n,$$ where $n$ is the number of observations. 
The BIC is one of the most-used criteria in the case of clustering \citep{fraley1998many, bouveyron2015discriminative, browne2011model, punzo2016parsimonious, schmutz2020clustering}. It is more parsimonious than the AIC since it generates greater penalties for the number of variables in the model. We opt for the BIC in this study.

\section{Numerical experiments}
\label{sec:num_exp}
This section aims to assess the performance of C-funHDDC on simulated data. A first set of experiments allows us to evaluate the efficiency of C-funHDDC itself and the influence of the hyper-parameters (type of initialization, intrinsic dimension, number of clusters, proportion of outliers). The second part of our experiments is dedicated to the comparison of C-funHDDC with competitors, namely {\it Functional Trimmed} \citep{garcia2005proposal} and FIF \citep{staerman2019functional}.

In all of our experiments, results are based on $100$ simulations. 
The C-funHDDC algorithm is implemented in \texttt{R}. For these simulated data sets, the computation time for one execution of the algorithm with {\it trimmed} initialization is between 30 seconds to 3 minutes on a 2.2 GHz Intel Xeon 4114 with 40 Go RAM.

\subsection{Simulated data}
\label{sec:simu}

Two simulated data sets are chosen, inspired by the work of \citep{preda2007regression, ferraty2003curves, bouveyron2015discriminative, schmutz2020clustering}. The two data sets have the same normal observations but have different abnormal observations. In the first one, \textit{dataset1}, all components of the multivariate curves are abnormal. In the second one, \textit{dataset2}, only  one component is abnormal. This latter case is most frequent in the industrial context. For each dataset, a sample of $n = 1005$ bivariate curves is simulated as follows: 
\begin{itemize}
    \item For both data sets, $1000$ normal curves are generated. The curves are composed of 250 curves derived from $K=4$ classes, according to the following model: for $t \in [1:21]$,
    $\text{ Class 1:} \;\; X_1(t) = U + (1-U)H_1(t) + \varepsilon(t),\; X_2(t) =  U + (0.5-U)H_1(t) +  \varepsilon(t), $
    $\text{ Class 2:} \;\; X_1(t) = U + (1-U)H_2(t) + \varepsilon(t), \;
    X_2(t) =  U + (0.5-U)H_2(t) +  \varepsilon(t),$
    $\text{ Class 3:} \;\; X_1(t) = U + (0.5-U)H_1(t) + \varepsilon(t), \;
    X_2(t) =  U + (1-U)H_2(t) +  \varepsilon(t),$
    $\text{ Class 4:} \;\; X_1(t) = U + (0.5-U)H_2(t) + \varepsilon(t), \;
    X_2(t) =  U + (1-U)H_1(t) +  \varepsilon(t),$ 
    where $U \sim \mathcal{U}{[0,1]}$ and $\varepsilon(t)$ is a Gaussian noise process. For all $t$, $\varepsilon(t)$ is Gaussian with zero mean and a variance equal to $0.25$, independent from $\epsilon(t')$ for $t'\neq t$. The functions $H_1$  and $H_2$ are the shifted triangular waveforms expressed as follows: for $t\in [1,21]$, $H_1(t) =  \max(6 -\mid t-7\mid, 0)$ and $H_2(t) =  \max(6 -\mid t-15\mid,0)$;
    \item For the first data set (where the two components are abnormal), we add $5$ abnormal curves partitioned into two groups. $3$ curves are simulated according to the following model: for $t \in [1:12]$, 
        $X_1(t) =  (0.5-U)H_1(t) + \sin(\pi t/2) + \varepsilon_2(t), \;
        X_2(t) =  (1-U)H_2(t)  + \sin(\pi t/2) + \varepsilon_2(t),$ 
    where $\varepsilon_2(t)$ is a Gaussian noise process independent from $\varepsilon(t)$ with zero mean and variance equal to $1$.
    2 curves are simulated from the following model: for $t \in [1,12]$, 
       $ X_1(t) = U + (1-U)H_3(t) + \varepsilon(t), \;
        X_2(t) = U + (0.5-U)H_3(t)  + \varepsilon(t)$
    where $H_3(t) = \max(6 -\mid t\textbf{1}_{\{ t<7\} } -7\mid,0)$.     We thus consider two types of outliers. The first type will be denoted in the following \emph{outlier1} and the second one \emph{outlier2}. We call this data set, {\it dataset1}.
    \item For the second data set (where one of the components is abnormal), we add to the normal data $5$ abnormal curves also partitioned into two groups. The first $3$ curves are defined as
        $X_1(t) = U + (0.5-U)H_1(t) + \varepsilon(t), \;
        X_2(t) = (1-U)H_2(t) + \sin(\pi t/2) + \varepsilon_2(t).$
     2 curves are simulated as 
        $X_1(t) = U + (1-U)H_3(t) + \varepsilon(t), \;
        X_2(t) =  U + (0.5-U)H_1(t) +  \varepsilon(t).$
    Similarly this data set contains two types of outliers. They will also be denoted respectively \emph{outlier1} and \emph{outlier2} in the following.
    This dataset is denoted {\it dataset2} below.
\end{itemize}

\begin{figure}[!ht]
    \centering
    \includegraphics[scale=0.35]{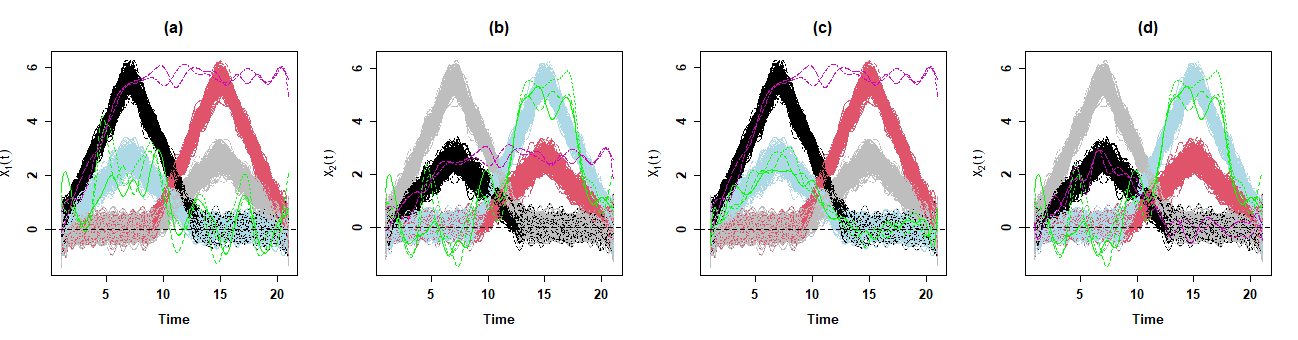}
    \caption{Simulated curves: {\it dataset1} ($2$ figures on the left) and {\it dataset2} ($2$ figures on the right). $X_1(t)$ is displayed on the left and $X_2(t)$ on the right. Each curve is colored according to its class: black for Class 1, red for Class 2, light-blue for Class 3 and gray for Class 4. Abnormal curves are colored respectively in green (outlier1) and in purple (outlier2).}
    \label{graph_data_simulated}
\end{figure}

The number of observation points is $T = 101$ with a regular frequency in $[1,21]$. To be consistent with \citep{bouveyron2015discriminative, schmutz2020clustering}, the reconstruction of the functional data is performed by cubic spline with $25$ basis functions. 
We represent the components of the two smoothed data in Figure~\ref{graph_data_simulated}. 

\subsection{Evaluation of C-funHDDC}
The goal of this section is to evaluate C-funHDDC and its ability to retrieve clusters and to detect outliers.
 
To evaluate the performance for both clustering and outlier detection, we use Adjusted Rand Index (ARI) \citep{rand1971objective}. Recall that ARI belong to $[0,1]$. An ARI equal to 1 indicates that the partition provided by the algorithm is perfectly aligned with the simulated one. Conversely, an ARI equal to 0 indicates that the two partitions are just some random matches.
We denote by $ARI_c$ the ARI index evaluating the quality of clustering for the normal curves. Similarly, $ARI_o$ evaluates the quality of outlier detection, computed on the whole set of normal and abnormal curves.

In the following, we focus on the evaluation of the initialization procedures and of the choice of the hyper-parameters. We first study the efficiency of the different types of initialization and the choice of the intrinsic dimension $d_k$. Next we evaluate the ability of the BIC to choose the number of clusters, and finally we assess the capacity of C-funHDDC to detect outliers.

\subsubsection{Initialization of the algorithm}
For the first experiment, the influence of different types of initialization, {\it trimmed, kmeans, random}, is evaluated. We consider both the ability of C-funHDDC to cluster the normal curves and to detect abnormal ones. A comparison of the convergence of the algorithms on one simulation is given in Supplementary Materials.
The number of clusters is fixed to $K = 4$. For the choice of $d_k$, the algorithm is run on a grid of positive integers (from $2$ to $10$), and the value of $d_k$ giving the best BIC is kept. Proceeding this way, we obtain $d_k = 2$ for each cluster. This choice of $d_k$ will be investigated in the next section.

\begin{figure}[!ht]
    \centering
    \includegraphics[scale=0.5]{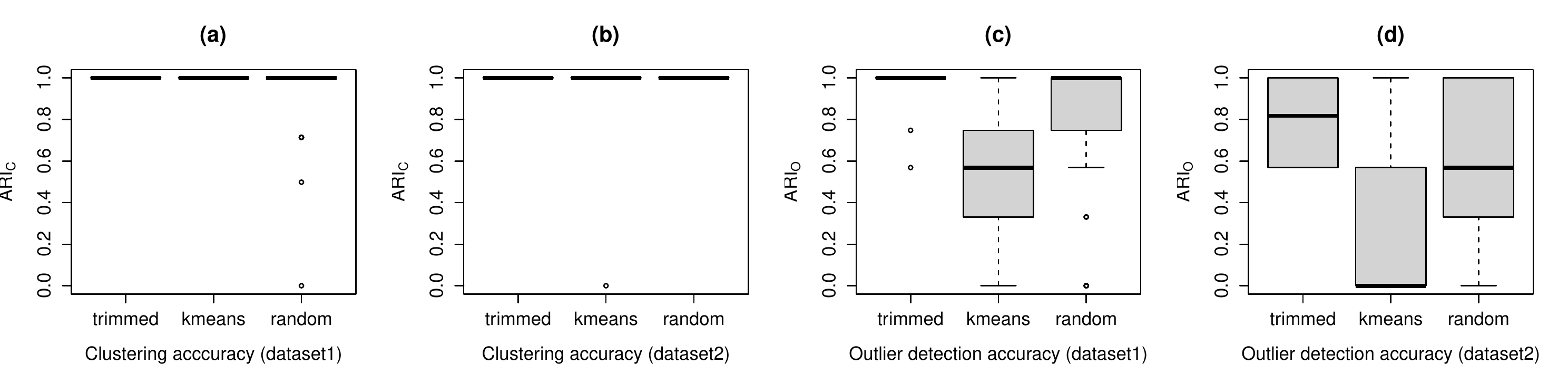}
    \caption{ARI partition with different types of initialization on $100$ simulations. The $2$ plots on the left give results for the clustering accuracy for {\it dataset1} (a) and {\it dataset2} (b). The plots on the right give results for the outliers detection accuracy for {\it dataset1} (c) and {\it dataset2} (d).}
    \label{graph_ari_boxplots}
\end{figure}

Figure \ref{graph_ari_boxplots} displays the boxplots of ARI using different types of initialization. It highlights that the {\it trimmed} approach is the best way to initialize the ECM algorithm for C-funHDDC. It perfectly groups the simulated normal data, and it clearly outperforms the other approaches in outlier detection. Consequently, {\it trimmed} will be used in the rest of the paper for initializing the ECM algorithm.

\subsubsection{Choice and influence of $d_k$ }
The intrinsic dimensions of each cluster are hyper-parameters of C-funHDDC which needs to be evaluated.
This can be done in two different ways: either by using a grid of values and selecting the best ones according to the BIC, or by using the Cattell scree test \citep{cattell1966scree}, with a threshold selected by BIC. The advantage of the second one is that it requires to select only one threshold for every cluster, whereas the grid search has a complexity that exponentially increases with $K$.

The number of clusters is fixed at $K=4$ and the ECM algorithm is initialized with {\it trimmed}.
The grid search explores all possible values of $(d_1,\ldots,d_K)$ in $\{2,10\}^K$.
For the Cattell threshold, all values from $0.1$ to $0.25$ with a step of $0.05$ are considered.
 
\begin{figure}[!ht]
    \centering
    \includegraphics[scale=0.5]{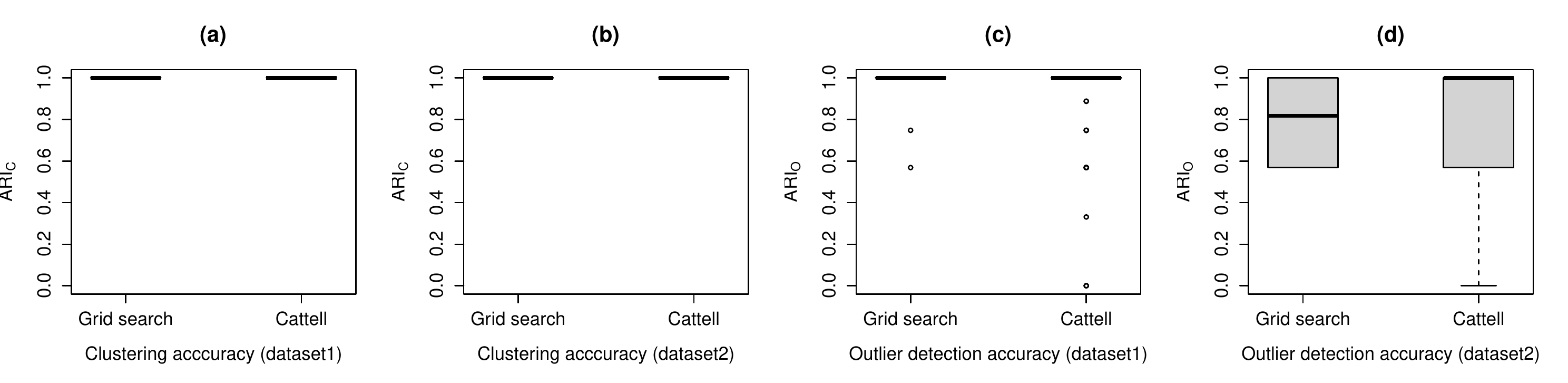}
    \caption{Clustering and outlier detection efficiency according to the manner of selecting the intrinsic dimension $d_k$:  Cattell scree test and grid search. Boxplots were obtained on $100$ simulations. The $2$ plots on the left give results for the clustering accuracy for {\it dataset1} (a) and {\it dataset2} (b). The plots on the right give results for the outliers detection accuracy for {\it dataset1} (c) and {\it dataset2} (d).}
    \label{graph_comparison_dk_cattell}
\end{figure}

Figure \ref{graph_comparison_dk_cattell} provides a comparison of the grid search and the Cattell scree test for selecting the best intrinsic dimensions. Both clustering and outliers detection performances are considered, through $ARI_c$ and $ARI_o$.

If the grid search gives more accurate results on {\it dataset1}, it is more nuanced with {\it dataset2} since the median of the results is better for Cattell. Nevertheless, the Cattell scree test sometimes obtains very low $ARI_o$, which may be a main drawback for an application to real data. Consequently, in the rest of the paper, the intrinsic dimensions $d_k$ will be selected by grid search. 

\subsubsection{Selection of the number of clusters}
In the following, we investigate the ability of the BIC to select the true number of clusters. Although the BIC has well-known consistence properties in the usual mixture model context, the presence of outliers may have a negative impact on its performance. BIC can have some difficulties for choosing the true number of clusters. We provide an example in Supplementary Materials, where funHDDC does not succeed in finding the right number of clusters.
The C-funHDDC algorithm is run for different values of $K$, $K \in \{1,2,3,4,5,6\}$. 
The intrinsic dimensions $(d_k)$ are selected using a grid search (on $\{2,3,4,5\}$). In each of the simulations, we obtained that $d_k=2$ for all $k$, whatever the final value of $K$.

Table~\ref{table:clusters_numbers1} reports the number of clusters chosen for $100$ simulations of both simulated datasets. Each line presents for each value of $k$, the number of times this value is chosen. Table~\ref{table:clusters_numbers1} shows clearly that the BIC tends to select 5 clusters rather than 4. An example of the clustering and outlier detection is displayed in Supplementary Materials.

\begin{table}[!ht]
\begin{minipage}{.45\linewidth}
\centering
    \begin{tabular}{ccc}
    \toprule
     Number of clusters  & {\it dataset1} & {\it dataset2} \\
     $(K)$ & & \\
    \midrule
    $1$ & $-$ & $-$\\ 
    $2$ & $1$ &  $-$ \\
    $3$ & $-$ &  $-$ \\
    $4$ & $9$  & $16$ \\
    $\mathbf{5}$ & $\mathbf{59}$  & $\mathbf{55}$ \\
    $6$ & $31$  & $29$ \\
    \bottomrule
    \end{tabular}
    \caption{ Best number of clusters selected by the BIC for $100$ simulations for each data set as a percentage, with choice of $d_k$ with grid search method. The highest proportions are displayed in bold type.}
    \label{table:clusters_numbers1}
\end{minipage}
\hspace{.5cm}
\begin{minipage}{.45\linewidth}
\centering
\begin{tabular}{cc}
    \toprule
    Number of clusters  & {\it Normal data}\\
    $(K)$ & \\
    \midrule 
    $1$ & $-$ \\
    $2$ & $1$ \\ 
    $3$ & $-$  \\
    $\mathbf{4}$ & $\mathbf{57}$  \\
    $5$ & $21$   \\
    $6$ & $21$  \\
    \bottomrule
    \end{tabular}
    \hspace{.5cm}
    \caption{ Best number of clusters selected by the BIC for $100$ simulations for normal data set as a percentage, with choice of $d_k$ with grid search method. The highest proportions are displayed in bold type.}
    \label{table:clusters_numbers2}
\end{minipage}

\end{table}

To explain the fact that 5 clusters are obtained with the BIC, we run a similar procedure on the dataset containing only normal curves of \textit{dataset1} and \textit{dataset2}, that is, removing all abnormal curves. Results are given in Table \ref{table:clusters_numbers2}. It highlights that the presence of 5 clusters in Table~\ref{table:clusters_numbers1} is due to the presence of abnormal data. Indeed if we look more precisely at the characteristics of the clusters, it appears that the algorithm tends to group the outliers into additional clusters. When abnormal data are discarded, the true number of cluster is chosen.

In conclusion, when the number of clusters is free, C-funHDDC tends to create additional clusters for the outliers. When the number of clusters is fixed to the number of clusters of normal curves, without an outlier, then the outliers are assigned to the closest cluster but well detected as outliers. The next section evaluates more precisely this ability to detect outliers when $K=4$.

\subsubsection{Outlier detection}
This section evaluates the ability of C-funHDDC to detect outliers, and the influence of different parameters: the proportion of outliers; and the level of heterogeneity among data. The number of clusters is fixed to $K = 4$ and the ECM algorithm is initialized with {\it trimmed}. In order to reduce computing time, the intrinsic dimensions are all fixed to $d_k=2$, which corresponds the values selected by grid search in the previous experiments.

\begin{figure}[!ht]
        \centering
        \includegraphics[scale=0.53]{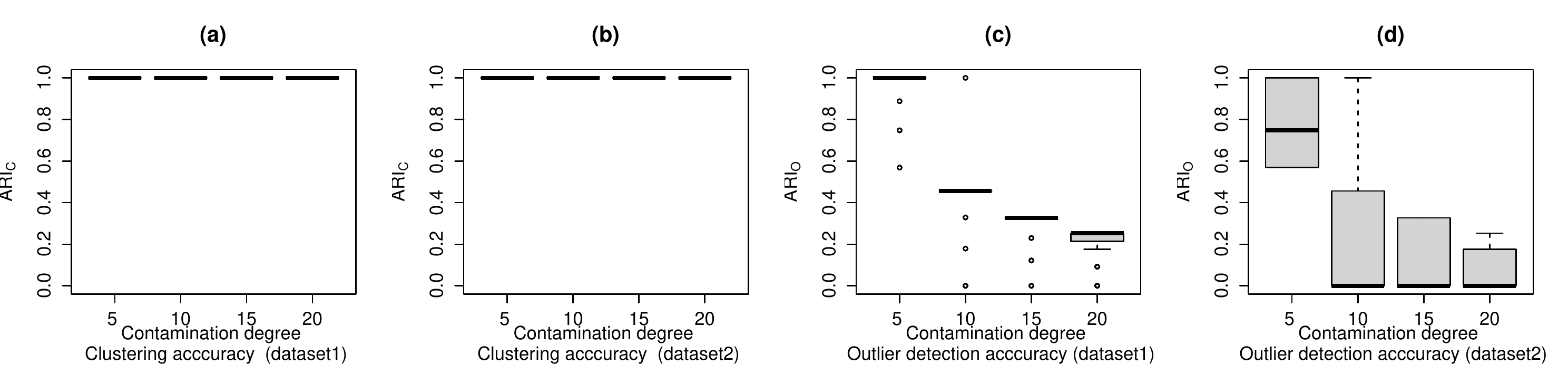}
        \caption{ARI for the clustering and for the outliers detection with respect to the degree of contamination (in per thousand) when $K=4$. The $2$ plots on the left give results for the clustering accuracy for {\it dataset1} (a) and {\it dataset2} (b). The $2$ plots on the right give results for the outliers detection accuracy for {\it dataset1} (c) and {\it dataset2} (d).}
    \label{graph_contamination_variation_2}
\end{figure}

\paragraph{Influence of the degree of contamination}
The initial \textit{dataset1} and \textit{dataset2} contain $5\mbox{\textperthousand}$ of outliers ($5$ abnormal curves for $1000$ normal ones). In this section, this proportion is gradually increased from $5\mbox{\textperthousand}$ to $20\mbox{\textperthousand}$ with step $5\mbox{\textperthousand}$, by adding the number of {\it outlier2} (initially set to 2). These outliers are not especially associated to one cluster, unlike the other type of outliers. Our results are given in Figure \ref{graph_contamination_variation_2}.

\begin{figure}[!ht]
    \centering
    \includegraphics[scale=0.5]{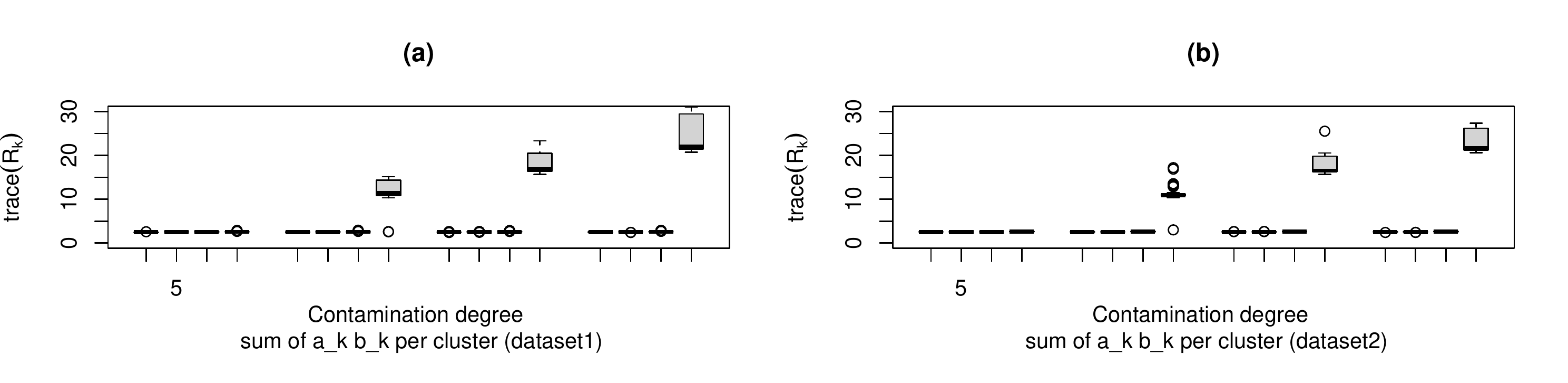}
    \caption{Values of $trace(R_k)$ for $1\leq k \leq K$ with respect to to the proportion of outliers, when $K=4$. Top plot (a) gives results for {\it dataset1} and bottom plot (b) for {\it dataset2}.}
    \label{graph_contamination_variance1}
\end{figure}

As we can see on Figures \ref{graph_contamination_variation_2}{\bf a} and \ref{graph_contamination_variation_2}{\bf c}, the clustering of normal data is not impacted by the variation of the number of outliers. Nevertheless, the ability of C-funHDDC algorithm to detect outliers decreases when their proportion increases (Figure \ref{graph_contamination_variation_2}{\bf b}, \ref{graph_contamination_variation_2}{\bf d}). 
The reason for this phenomenon is that when the proportion of outliers becomes too large, for $K$ fixed to 4, the model prefers to increase the variance of the clusters and considers the abnormal curves as normal rather than as an outlier. This is illustrated by Figure~\ref{graph_contamination_variance1}. Recall that the variance structure in each cluster $k$ is given by the matrix $R_k$ defined in \eqref{eq:HDDC}. Figure~\ref{graph_contamination_variance1} plots the value of the trace of the matrices $\{R_k, \;k=1,\dots,K\}$ with respect to the proportion of outliers. Indeed, one of the clusters has a variance which increases with the proportion of outliers.

\begin{figure}[!ht]
    \centering
    \includegraphics[scale=0.35]{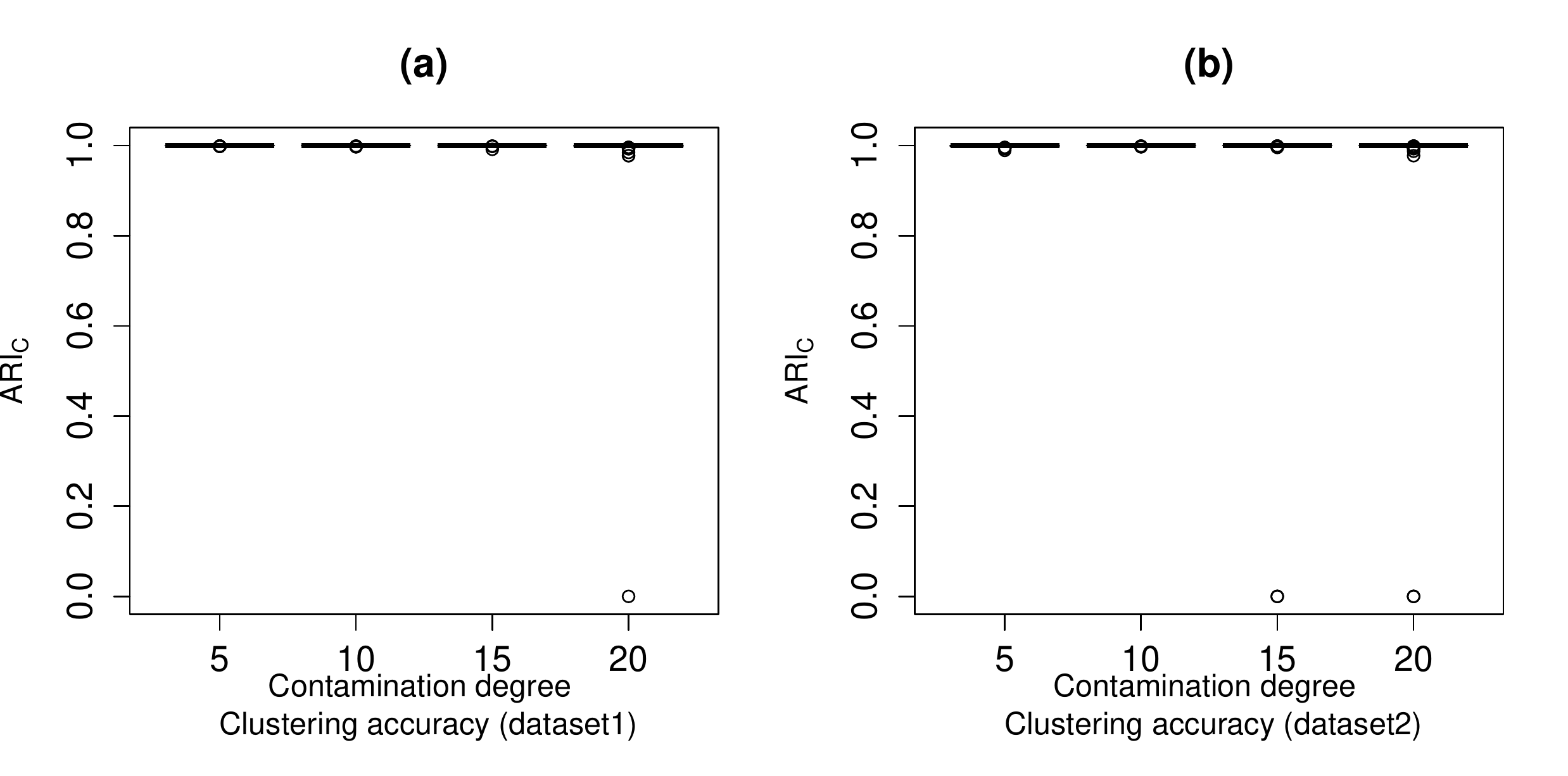}
    \caption{ARI for the clustering with respect to the degree of contamination (per thousand) when $K=5$. Plot (a) gives results for {\it dataset1} and plot (b) for {\it dataset2}. As outliers are associated to a cluster, the ARI index measures here both clustering accuracy and outlier detection accuracy.}
    \label{graph_sample_variation_boxplots_K=5}
\end{figure}

\begin{figure}[!ht]
    \centering
    \includegraphics[scale=0.55]{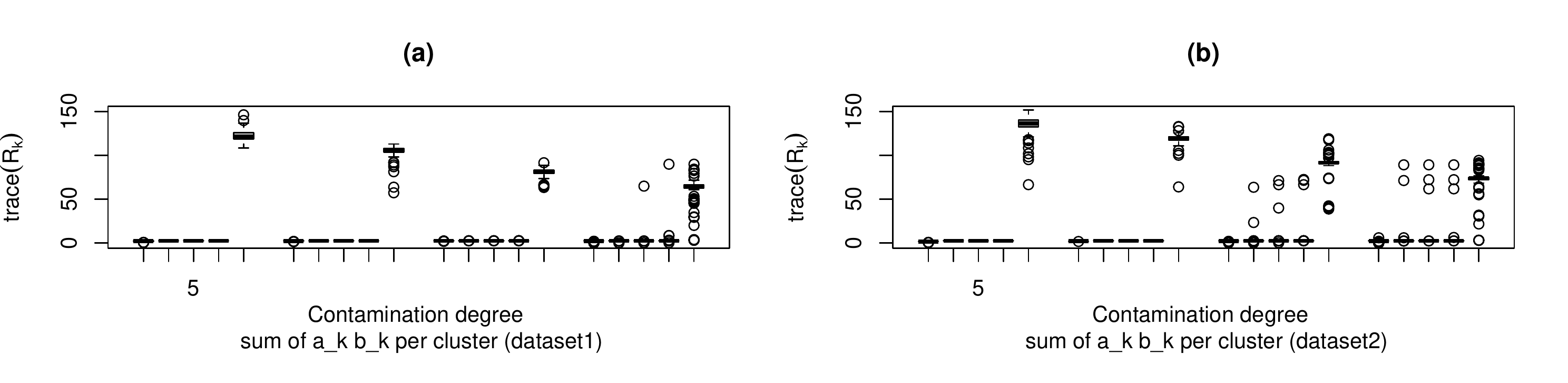}
    \caption{Values of $trace(R_k)$ for $1\leq k \leq K$ with respect to to the proportion of outliers, when $K=5$. Top plot (a) gives results for {\it dataset1} and bottom plot (b) for {\it dataset2}.}
    \label{graph_contamination_variance2}
\end{figure}

Now consider the case where the number of clusters is free. Then the BIC selects $K=5$ and all the outliers are grouped into an additional cluster. This is illustrated in Figure~\ref{graph_sample_variation_boxplots_K=5} which provides the ARI results for $K=5$. All the abnormal data are perfectly grouped in the same cluster. Similarly, we can observe the evolution of the variance within clusters. Figure~\ref{graph_contamination_variance2}  displays the value of the trace of the covariance matrices $\{R_k, \;k=1,\dots,K\}$ with respect to the proportion of outliers. In this setting, the variances of the 4 clusters containing normal curves remain stable. By contrast, the variance of the cluster containing the outliers is much higher. This is due to the presence of two types of outliers ({\it outlier1} and {\it outlier2}).

\begin{figure}[!ht]
    \centering
    \includegraphics[scale=0.5]{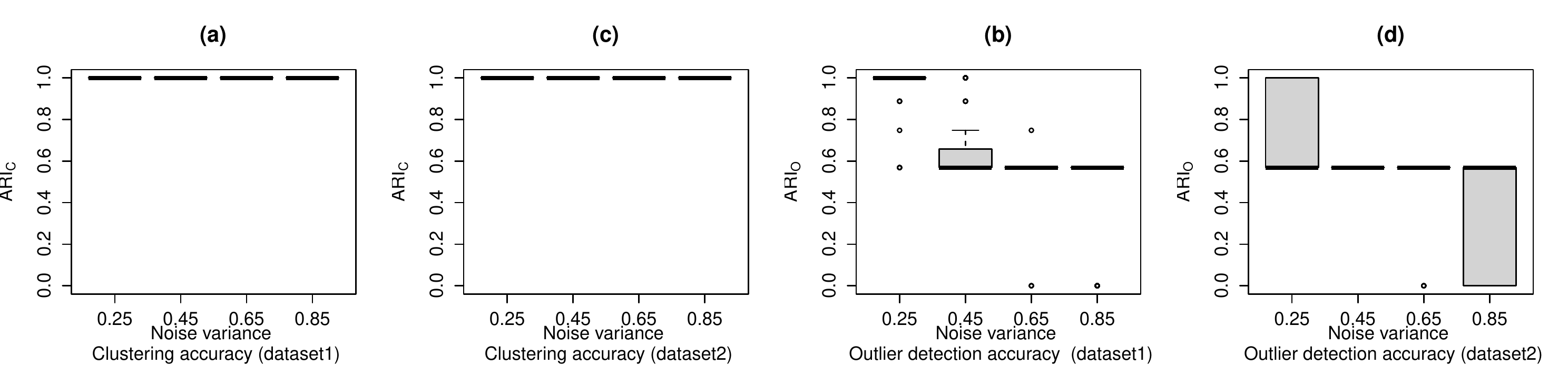}
    \centering
    \caption{ARI for the clustering and for the outliers detection with respect to the variance within clusters. The $2$ plots on the left give results for the clustering accuracy for {\it dataset1} (a) and {\it dataset2} (b). The plots on the right give results for the outliers detection accuracy for {\it dataset1} (c) and {\it dataset2} (d).}
    \label{graph_noise_variation_boxplots}
\end{figure}

\paragraph{Influence of the level of heterogeneity within clusters}
Since the clusters are relatively well separated in \textit{dataset1} and \textit{dataset2}, we now investigate the influence of the degree to which they are mixed. The noise variance within each cluster is increased. The variance of variable $\varepsilon(t)$ defined in Section~\ref{sec:simu} is increased from 0.25 to 0.85 with a step of $0.2$. For each variance value, $100$ samples are simulated and C-funHDDC is performed with the same setting as selected in the previous experiments. $K$ is fixed at $4$. Plots of $ARI_c$ and $ARI_o$ with respect to the noise variance are displayed in Figure \ref{graph_noise_variation_boxplots}. The change of the noise within variances does not affect the clustering of the normal data. By contrast, the ability to detect outliers decreases when the cluster variance increases.

\subsection{Comparison with competitors}
In this section, C-funHDDC is compared with some competitors according to the ability of detecting outliers. The competitors are {\it Functional Trimmed} (FT) \citep{garcia2005proposal} and {\it Functional Isolation Forest} (FIF) \citep{staerman2019functional}. 

The FT method is defined as the kmeans method applied on the basis expansion coefficients of the curves, after having removed of proportion $\alpha$ of extreme observations. This hyper-parameter $\alpha$ has to be set, and there is no automatic way to do that in the present unsupervised context. Consequently, in the experiments, we use different values of $\alpha = {0.001, 0.005, 0.01}$, including the true one ($0.005$).

The FIF algorithm is associated with a score that gives the abnormality of each piece of data. This score depends on a {\it path length} which reflects the number of cuts necessary to isolate a datum. We choose to use $L_2$ scalar product of derivatives and to project the observations on a Gaussian-wavelets dictionary. FIF algorithm uses a threshold to specify whether the data are normal. To define this threshold, the authors introduced the same hyper-parameter $\alpha$ (named \emph{contamination} in the paper). That is, $\alpha$ is an estimation of the proportion of outliers in the data set. Here also, we choose $3$ values for this contamination parameter $\{0.001, 0.005, 0.01\}$.

These algorithms are already implemented by the authors. The FT algorithm is implemented in \texttt{R}. It is available on request from the authors. The FIF algorithm is implemented in \texttt{Python} and the code can be accessed at: \url{https://github.com/Gstaerman/FIF}.

Note that these two approaches require the knowledge of the proportion of abnormal curves in the dataset to set a “threshold of normality”. This is not the case for C-funHDDC, which does not require to know the estimation of the proportion of outliers before declaring that a datum is normal or abnormal. 
In real data applications, it is not possible to independently estimate the proportion of outliers. Therefore, the adaptivity of C-funHDDC procedure is a major advantage with respect to FT and FIF.

\begin{figure}[!ht]
    \centering
    \includegraphics[scale=0.5]{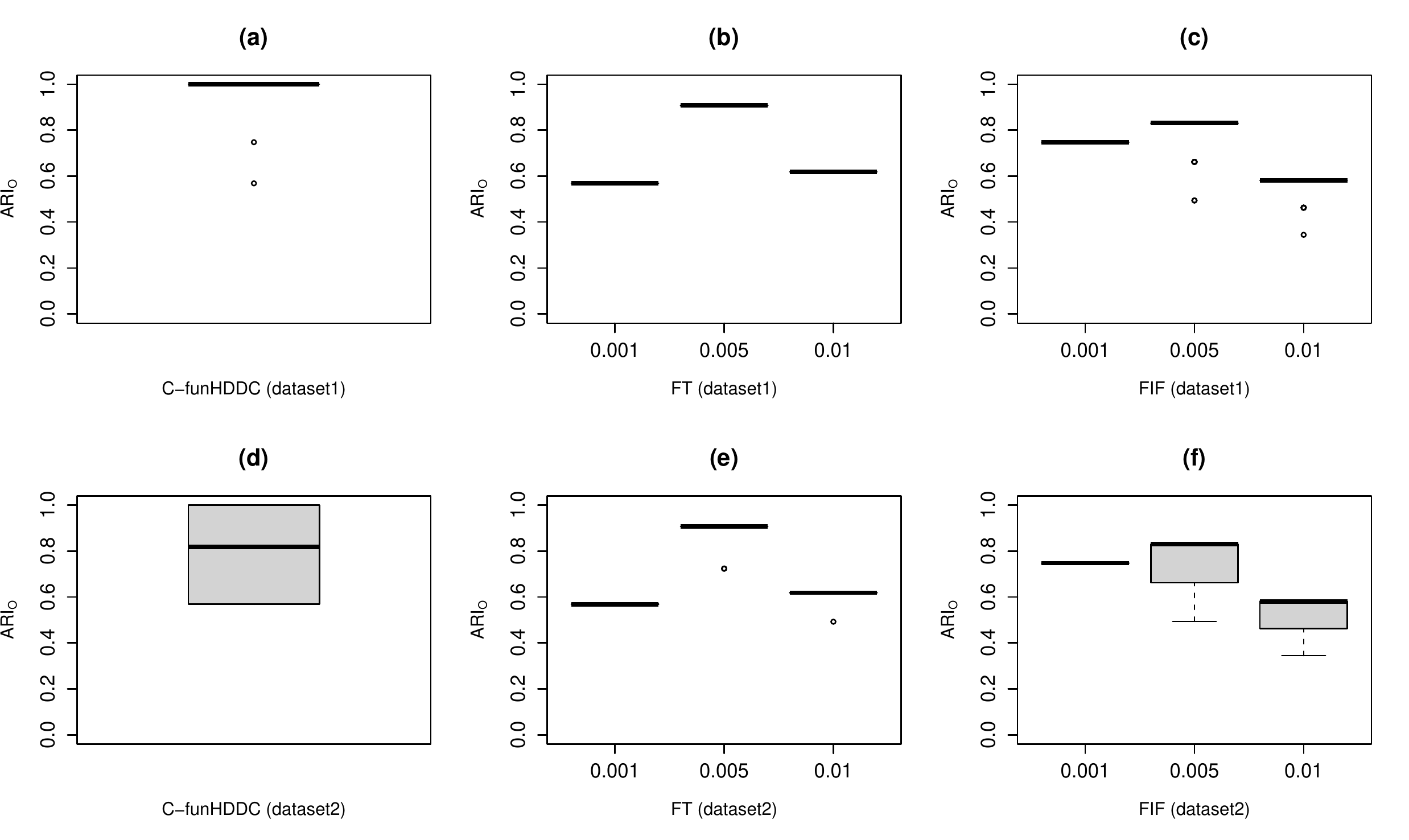}
    \caption{Outlier detection accuracy ($ARI_o$) for C-funHDDC, FT and FIF. Top plots give results for {\it dataset1}: (a) C-funHDDC (b) FT  (c) FIF. Bottom plots give results for {\it dataset2}: (d) C-funHDDC (e) FT  (f) FIF.}
    \label{graph_fif_boxplots}
\end{figure}

Figure \ref{graph_fif_boxplots} presents the results for the three methods.
In the light of these results, the best performances are obtained with C-funHDDC. FIF and FT perform both correctly when the proportion of outliers is known, but the results are worse when this proportion is wrongly fixed.

\section{Application to real data}
\label{sec:real_data}

\subsection{Data and preprocessing}
The data under analysis in this section consist of 569 measurements by four sensors. 
Figure \ref{graph_data2} represents a sample of 148 measurements (one plot per sensor).
The acquisition frequency is 100Hz.
The time duration of each measurement is not constant. The number of observed points for each measurement varies from 2199 to 10675, with a median equal to 5144. Consequently, a preliminary time normalization is performed so as to bring all the time periods to $[0,1]$. This is done by rescaling the time as follows: $t = \frac{t - t_{min}}{t_{max} - t_{min}}$, where $t_{min}$ and $t_{max}$ are the starting and ending time of each measurement. The basis functions and the number of basis functions are chosen empirically. The functional data are reconstructed using a cubic spline basis with 25 basis functions. In this data set, 45 outliers have been defined by an expert, and our goal is to correctly detect them without detecting too many false positives.

\subsection{Results}
C-funHDDC is applied to the industrial data set by selecting the number $K$ of clusters and the intrinsic dimension $d_k$ by BIC.
The range of explored values is $\{2, 3, 4\}$ for $K$ and it is $\{2, 3 ,\dots, 10\}$ for $d_k$. 
The range of values is limited to 10 for $d_k$ in order to favor the selection of parsimonious models.
The best number of clusters selected by the BIC is $K = 3$ with $(d_1, d_2, d_3) = (10, 10, 6)$.

As described in Section \ref{sec:num_exp}, when the outliers are numerous, C-funHDDC gathers them into the same cluster. This cluster becomes a cluster of outliers. This is what happens in this industrial application. Indeed, no datum is detected as  abnormal observation, but one of the clusters contains all of the outliers. The results are presented in Table~\ref{table:conf_matrix1} as a confusion matrix. An additional description of the clusters is provided in Supplementary Materials.

\begin{table}[!ht]
 \centering
\begin{tabular}{@{}cc cc@{}}
\multicolumn{1}{c}{} &\multicolumn{1}{c}{} &\multicolumn{2}{c}{Predicted} \\ 
\cmidrule(lr){3-4}
\multicolumn{1}{c}{} & 
\multicolumn{1}{c}{} & 
\multicolumn{1}{c}{normal} & 
\multicolumn{1}{c}{outlier} \\ 
\cline{2-4}
\multirow[c]{2}{*}{\rotatebox{90}{~Reference}}
& normal  & 522 & 2   \\[1.5ex]
& outlier  & 0   & 45 \\ 
\cline{2-4}
\end{tabular}
\vspace{0.15cm}
\caption{Confusion matrix using C-funHDDC model}
\label{table:conf_matrix1}
\end{table}


C-funHDDC detects all outliers, with only $2$ false positives. For fixed values of $K$ and $d_k$, the computation time for one execution of C-funHDDC is $3.72$ minutes.

In order to evaluate the performance of C-funHDDC, Table \ref{table:conf_matrix2} and \ref{table:conf_matrix3} present the results obtained respectively with FT and with FIF.
The problem with FT and FIF is that they need to fix the proportion of outliers, which is clearly unknown in real application. In this application three proportions are used ($ \alpha \in \{0.01, 0.08, 0.15 \}$), among them the true one ($0.08$), which is advantageous for FT and FIF. Nevertheless, even with the true proportion of outliers, the results are relatively poor with FT and FIF, in particular when comparing these with the results of C-funHDDC. The computation times for one execution of the FIF and FT algorithm are respectively $08.24$ minutes and $19.08$ seconds.

\begin{table}[!ht]
\centering
\begin{tabular}{ccc}
\begin{tabular}{@{}cc cc@{}}
\multicolumn{4}{c}{$\alpha=0.01$}\\
\multicolumn{1}{c}{} &\multicolumn{1}{c}{} &\multicolumn{2}{c}{Predicted} \\ 
\cmidrule(lr){3-4}
\multicolumn{1}{c}{} & 
\multicolumn{1}{c}{} & 
\multicolumn{1}{c}{normal} & 
\multicolumn{1}{c}{outlier} \\ 
\cline{2-4}
\multirow[c]{2}{*}{\rotatebox{90}{Reference}}
& normal  & 524 & 0   \\[1.5ex]
& outlier  & 39   & 6 \\ 
\cline{2-4}
\end{tabular}
&
\begin{tabular}{@{}c cc@{}}
\multicolumn{3}{c}{$\alpha=0.08$}\\
\multicolumn{1}{c}{} &\multicolumn{2}{c}{Predicted} \\ 
\cmidrule(lr){2-3}
\multicolumn{1}{c}{} & 
\multicolumn{1}{c}{normal} & 
\multicolumn{1}{c}{outlier} \\ 
\cline{1-3}
 normal  & 505 & 19   \\[1.5ex]
 outlier  & 18   & 27 \\ 
\cline{1-3}
\end{tabular}

&

\begin{tabular}{@{}c cc@{}}
\multicolumn{3}{c}{$\alpha=0.15$}\\
\multicolumn{1}{c}{} &\multicolumn{2}{c}{Predicted} \\ 
\cmidrule(lr){2-3}
\multicolumn{1}{c}{} & 
\multicolumn{1}{c}{normal} & 
\multicolumn{1}{c}{outlier} \\ 
\cline{1-3}
 normal  & 461 & 55   \\[1.5ex]
 outlier  & 14   & 31 \\ 
\cline{1-3}
\end{tabular}

\end{tabular}
\caption{Confusion matrix with FT}
\label{table:conf_matrix2}
\end{table}

\begin{table}[!ht]
 \centering
 \begin{tabular}{ccc}
\begin{tabular}{@{}cc cc@{}}
\multicolumn{4}{c}{$\alpha=0.01$}\\
\multicolumn{1}{c}{} &\multicolumn{1}{c}{} &\multicolumn{2}{c}{Predicted} \\ 
\cmidrule(lr){3-4}
\multicolumn{1}{c}{} & 
\multicolumn{1}{c}{} & 
\multicolumn{1}{c}{normal} & 
\multicolumn{1}{c}{outlier} \\ 
\cline{2-4}
\multirow[c]{2}{*}{\rotatebox{90}{Reference}}
& normal  & 521 & 3  \\[1.5ex]
& outlier  & 41   & 4 \\ 
\cline{2-4}
\end{tabular}

&

\begin{tabular}{@{}c cc@{}}
\multicolumn{3}{c}{$\alpha=0.08$}\\
\multicolumn{1}{c}{} &\multicolumn{2}{c}{Predicted} \\ 
\cmidrule(lr){2-3}
\multicolumn{1}{c}{} & 
\multicolumn{1}{c}{normal} & 
\multicolumn{1}{c}{outlier} \\ 
\cline{1-3}
 normal  & 491 & 33   \\[1.5ex]
 outlier  & 31   & 14 \\ 
\cline{1-3}
\end{tabular}

&

\begin{tabular}{@{}c cc@{}}
\multicolumn{3}{c}{$\alpha=0.15$}\\
\multicolumn{1}{c}{} &\multicolumn{2}{c}{Predicted} \\ 
\cmidrule(lr){2-3}
\multicolumn{1}{c}{} & 
\multicolumn{1}{c}{normal} & 
\multicolumn{1}{c}{outlier} \\ 
\cline{1-3}
normal  & 460 & 64   \\[1.5ex]
outlier  & 22   & 23 \\ 
\cline{1-3}
\end{tabular}

\end{tabular}
\caption{Confusion matrix with FIF}
\label{table:conf_matrix3}
\end{table}

\section{Conclusion}

In this study, we are interested in clustering and in detecting outliers in multivariate functional data. We propose a contaminated Gaussian mixture model, derived from a decomposition of the data in functional basis. The model is based on a functional latent block model. For each cluster, a parameter controlling the proportion of outliers and one specifying the variance inflation factor from normal data are introduced to take into account the presence of outliers. The model was initially proposed in \cite{punzo2016contaminatedmixt} and we have extended it to a multivariate functional framework following \cite{jacques2013funclust}. An  ECM algorithm enables us to infer the parameters of the model.  Numerical experiments have shown that C-funHDDC has a very satisfactory behavior for clustering and detecting outliers in multivariate functional data. In particular, to our knowledge, it is the only adaptive procedure with respect to the proportion of outliers for such data. 

In the real data application, the outliers form a single cluster. Using a robust approach makes this possible, while non-robust approaches cannot identify the outliers. We do not detect any outliers associated with a cluster, but it is likely to happen later in the application. For the quality assessment of measurements from sensors, the C-funHDDC model makes it possible to decide whether a datum represented by several curves is correct or not.

A perspective gained by this study suggests the benefits of extending the model to non-Gaussian distributions. Indeed, if the Gaussian assumption is often realistic, it could be interesting to extend this work to other types of distribution, inspired for instance by \citep{forbes2014new}, or even by \citep{morris2019asymmetric} and \citep{punzo2021multiple} in the contaminated case.

Moreover, this work is motivated by an industrial application, where an online procedure is needed. We therefore plan to extend this modeling using scalable clustering and online mixture, as in \citep{bellas2013model}.

\section*{Acknowledgment}
We thank Arpege Master K for their financial support and for their help in analyzing the outliers detection from an expert point of view. We are also grateful to the reviewers for their insightful comments to improve this paper.

\appendix

\section{Proofs}
\label{sec:sample:appendix}



Suppose that the parameters of the model are estimated at the iteration $m$ and denote $\theta^{(m)}$ the corresponding estimation. 
Recall that our goal is to maximize $Q(\theta\mid \theta^{(m)})$, where $Q(\theta\mid \theta^{(m)}) = \E\left(\ell_c(\textbf{c},\textbf{z},\textbf{v},\theta) \mid \textbf{c}, \theta^{(m)} \right)$.  


\subsection{First Maximization step (CM1)}

At this step, we fix $\theta_2$ at $\theta_2^{(m)}$ and we compute $\theta^{(m+1)}$ as the value of $\theta_1$  that maximizes $Q(\theta\mid\theta^{(m)})$.

\paragraph{Estimation of $\pi_k$} 
Denote $Q_1(\theta\mid \theta^{(m)}) = \E \left(\ell_{c_1} (\textbf{c},\textbf{z}) \mid \theta^{(m)} \right) = \sum_{i=1}^n \sum_{k=1}^K \E( z_{ik}\mid c_i,\theta^{(m)}) \log(\pi_k) = \sum_{i=1}^n \sum_{k=1}^K z_{ik}^{(m)} \log(\pi_k) $. 
In order to take into account the relation which links $\{\pi_k, k=1,\ldots,K\}$, we introduce the Lagrange multiplier by using the function $Q_{L1} (\theta,\lambda\mid \theta^{(m)}) = Q_1 (\theta\mid \theta^{(m)}) - \lambda\left( \sum_{k=1}^K \pi_k -1 \right)$. The partial derivative with respect to $\pi_k$ and $\lambda$ are: 
\begin{align*}
    \frac{\partial Q_{L1}}{\partial \pi_k} = \sum_{i=1}^n \frac{t_{ik}^{(m)}}{\pi_k} - \lambda 
    \qquad \text{ and } \qquad 
    \frac{\partial Q_{L1}}{\partial \lambda} =  \sum_{k=1}^K \pi_k -1\,.
\end{align*}
By canceling these expressions, we get:
\begin{align}
  \pi_k^{(m+1)} = \frac{\sum_{i=1}^n t_{ik}^{(m)}}{n}\,.
\end{align}

\paragraph{Estimation of $\beta_k$}
Let 
\begin{align*}
    Q_2(\theta\mid \theta^{(m)}) = \E \left(\ell_{c_2} (\textbf{c},\textbf{z},\textbf{v},\theta) \mid \bf{c}, \theta^{(m)} \right)\
    = \sum_{i=1}^n \sum_{k=1}^K \left[ t_{ik}^{(m)} s_{ik}^{(m)}\log(\beta_k) +  t_{ik}^{(m)} (1 - s_{ik}^{(m)}) \log(1-\beta_k) \right]\,.
\end{align*}
By differentiating $Q_2(\theta\mid \theta^{(l)})$ with respect to $\beta_k$, we have:
\begin{align}
    \nonumber \frac{Q_2(\theta\mid \theta^{(m)})}{\partial\beta_k} = 0 & \iff \sum_{i=1}^n \left[ \frac{t_{ik}^{(m)} s_{ik}^{(m)}}{\beta_k} - \frac{t_{ik}^{(m)} (1 - s_{ik}^{(m)}}{1- \beta_k} \right] = 0\\
    &\iff \beta_k^{(m+1)} = \frac{\sum_{i=1}^n t_{ik}^{(m)} s_{ik}^{(m)}}{\sum_{i=1}^n t_{ik}^{(m)}}\,.
\end{align}

\paragraph{ Estimation of $\mu_k$}
Since $g(c_i,\mu_k,\Sigma_k) = \frac{1}{(2\pi)^{\frac{B}{2}} |\Sigma_k|^\frac{1}{2}} \exp{ \{-\frac{1}{2}(c_i - \mu_k)^\prime\Sigma_k^{-1}(c_i -\mu_k)\}},$
\begin{align*}
    \ell_{c_3}(\textbf{c},\textbf{z},\textbf{v},\theta )  &= -\frac{1}{2}\sum_{k=1}^K \sum_{i=1}^n \Bigl[  B  z_{ik}(1-v_{ik})  \log \eta_k + z_{ik}  \log |\Sigma_k|  \\ & \qquad \qquad\qquad \quad+   z_{ik} \left( v_{ik} + \frac{1-v_{ik}}{\eta_k} \right) (c_i - \mu_k)^\prime\Sigma_k^{-1}(c_i-\mu_k)   \Bigr]  - \frac{n}{2} B\log 2\pi\,.
\end{align*}
The derivative with respect to $\mu_k$ is:
\begin{align*}
    \frac{\partial \ell_{c_3}(\textbf{c},\textbf{z},\textbf{v},\theta )}{\partial\mu_k} =& \sum_{i=1}^n z_{ik}\left( v_{ik} + \frac{1-v_{ik}}{\eta_k} \right) \frac{\partial}{\partial\mu_k} \left[ (c_i -\mu_k)^\prime \Sigma_k^{-1} (c_i-\mu_k)  \right] \\
    =& \sum_{i=1}^n z_{ik}\left( v_{ik} + \frac{1-v_{ik}}{\eta_k} \right) 2 (c_i -\mu_k) \Sigma_k^{-1} \,.
\end{align*}
Consequently,
\begin{align*}
    \frac{\partial \ell_{c_3}(\textbf{c},\textbf{z},\textbf{v},\theta )}{\partial\mu_k} = 0 & \iff \mu_k =  \frac{1}{\sum_{i=1}^n z_{ik}\left( v_{ik} + \frac{1-v_{ik}}{\eta_k} \right)} \sum_{i=1}^n z_{ik}\left( v_{ik} + \frac{1-v_{ik}}{\eta_k} \right) c_i\,.
\end{align*}
Thereby,
\begin{equation}
    \mu_k^{(m+1)} =  \frac{1}{\sum_{i=1}^n t_{ik}^{(m)}\left( s_{ik}^{(m)} + \frac{1-s_{ik}^{(m)}}{\eta_k^{(m)}} \right)} \sum_{i=1}^n t_{ik}^{(m)}\left( s_{ik}^{(m)} + \frac{1-s_{ik}^{(m)}}{\eta_k^{(m)}} \right) c_i\,.
\end{equation}

\textbf{Estimation of $Q_k$ and the parameters $a_{kl}$ and $b_k$ }

As in the previous section, let $n_k$ be the number of curves which belong to the cluster $k$, $n_k =  \sum_{i=1}^n z_{ik} $. Let also $n_{k_g}$ be the number of  normal data in the cluster $k$. That is, $n_{k_g} =  \sum_{i=1}^n z_{ik} v_{ik}  $. The number of  abnormal data in the cluster $k$ is $n_k -n_{k_g} =  \sum_{i=1}^n z_{ik} (1-v_{ik})$.

Recall that the covariance matrix $\Sigma_k$ of the $k^\text{th}$ cluster is defined by: \[\Sigma_k = W^{-\frac{1}{2}} Q_k R_k Q_k^\prime W^{-\frac{1}{2}}.\] Since the matrices $R_k, Q_k$ and $W^{-\frac{1}{2}}$ are invertible, we can easily write  \[\Sigma_k^{-1} = W^{\frac{1}{2}} Q_k R_k^{-1} Q_k^\prime W^{\frac{1}{2}}.\]

As $|\Sigma_k | = | W^{-\frac{1}{2}}| |Q_k| |R_k| |Q_k^\prime| |W^{-\frac{1}{2}}| = |R_k| |Q_k Q_k^\prime| | W|^{-1} = |R_k| = \prod_{l=1}^{d_k} a_{kl} \prod_{l=d_k +1}^B b_k $, we get:
\begin{multline}
   \label{eqn:lc3}    \ell_{c_3}(\textbf{c},\textbf{z},\textbf{v},\theta ) = -\frac{1}{2}\sum_{k=1}^K \sum_{i=1}^n 
    \left[  B (n_k - n_{k_g}) \log \eta_k + n_k \log \left( \prod_{l=1}^{d_k} a_{kl} \prod_{l=d_k +1}^B b_k \right) \right.\\ 
  \left. + z_{ik} \left( v_{ik} + \frac{1-v_{ik}}{\eta_k} \right)  (c_i - \mu_k)^\prime W^{\frac{1}{2}} Q_k R_k^{-1} Q_k^\prime W^{\frac{1}{2}}(c_i-\mu_k) \right] - n B\log (2\pi)\,.
\end{multline}

Since $(c_i - \mu_k)^\prime W^{\frac{1}{2}} Q_k R_k^{-1} Q_k^\prime W^{\frac{1}{2}}(c_i-\mu_k$) is a  scalar,
\begin{align*}
    (c_i - \mu_k)^\prime W^{\frac{1}{2}} Q_k R_k^{-1} Q_k^\prime W^{\frac{1}{2}}(c_i-\mu_k) &= \tr\left((c_i - \mu_k)^\prime W^{\frac{1}{2}} Q_k R_k^{-1} Q_k^\prime W^{\frac{1}{2}}(c_i-\mu_k)\right) \\ 
    &= \tr\left[ \left((c_i - \mu_k)^\prime W^{\frac{1}{2}} Q_k\right) \left( R_k^{-1} Q_k^\prime W^{\frac{1}{2}}(c_i-\mu_k)\right)\right] \\
    &= \tr\left[ \left( R_k^{-1} Q_k^\prime W^{\frac{1}{2}}(c_i-\mu_k)\right) \left((c_i - \mu_k)^\prime W^{\frac{1}{2}} Q_k\right) \right]\,.
\end{align*}

Plugging this expression into \eqref{eqn:lc3} leads to:
\begin{multline*}
     \sum_{k=1}^K \left[  \sum_{i=1}^n z_{ik} \left( v_{ik} + \frac{1-v_{ik}}{\eta_k} \right)  \left( (c_i - \mu_k)^\prime W^{\frac{1}{2}} Q_k R_k^{-1} Q_k^\prime W^{\frac{1}{2}}(c_i-\mu_k) \right) \right]   \\
     = \sum_{k=1}^K  n_k  \, \tr\left[ \left( R_k^{-1} Q_k^\prime W^{\frac{1}{2}}\right) \Bigl(\frac{1}{n_k}   \sum_{i=1}^n z_{ik} \bigl( v_{ik} + \frac{1-v_{ik}}{\eta_k} \bigr) (c_i-\mu_k)^\prime (c_i - \mu_k) \Bigr) \left( W^{\frac{1}{2}} Q_k\right) \right]\,.
\end{multline*}
Denote \[
    H_k = \frac{1}{n_k} \left(  \sum_{i=1}^n z_{ik} \left( v_{ik} + \frac{1-v_{ik}}{\eta_k} \right) (c_i-\mu_k)^\prime (c_i - \mu_k) \right)\,.
\]
Hence $\ell_{c_3}$ can be written as:
\begin{align*}
    \ell_{c_3}(\textbf{c},\textbf{z},\textbf{v},\theta ) & =  -\frac{1}{2}\sum_{k=1}^K  
    \Bigl\{  B (n_k - n_{k_g}) \log \eta_k + n_k  \left( \sum_{l=1}^{d_k} \log (a_{kl}) + \sum_{l=d_k +1}^B \log(b_k) \right) \\
    & \qquad \qquad  + n_k tr\left(  R_k^{-1} Q_k^\prime W^{\frac{1}{2}} H_k W^{\frac{1}{2}} Q_k\right)  \Bigr\} - \frac{nB}{2}\log (2\pi) \\
    & = -\frac{1}{2}\sum_{k=1}^K  \Bigl\{ B (n_k - n_{k_g}) \log \eta_k + n_k  \left( \sum_{l=1}^{d_k} \log (a_{kl}) + \sum_{l=d_k +1}^B \log(b_k) \right)\\ 
    &\qquad\qquad + n_k tr\left[ \sum_{l=1}^{d_k} \frac{q_{kl}^\prime W^{\frac{1}{2}} H_k W^{\frac{1}{2}} q_{kl}}{a_{kl}} + \sum_{l=d_k+1}^B \frac{q_{kl}^\prime W^{\frac{1}{2}} H_k W^{\frac{1}{2}} q_{kl}}{b_k}\right]  \Bigr\} -\frac{nB}{2}\log(2\pi)\,.
\end{align*}

Our objective is to maximize $\ell_{c_3}$ under the constraint $q_{kl}^\prime q_{ku} = 0$ if $l \neq u$. It is equivalent to look for a saddle point of the Lagrange function:
\begin{align}\label{eq_lagrange_for_Qk}
    \mathcal{L}(\textbf{c},\textbf{z},\textbf{v},q_{kl}) = -2\ell_{c_3}(\textbf{c},\textbf{z},\textbf{v},\theta ) - \sum_{l=1}^B w_{kl} \left( q_{kl}^\prime q_{kl} -1 \right)\,.
\end{align}
If we differentiate with respect to $q_{kl}$, we obtain:
\begin{align*}
    \nabla_{q_{kl}} \mathcal{L}(\textbf{c},\textbf{z},\textbf{v},q_{kl}) &= 2 n_k \left[ \frac{W^{\frac{1}{2}} H_k W^{\frac{1}{2}} q_{kl}}{a_{kl}} + \frac{W^{\frac{1}{2}} H_k W^{\frac{1}{2}} q_{kl}}{b_k}\right] - 2 w_{kl}q_{kl} \\
    &= 2 n_k \frac{W^{\frac{1}{2}} H_k W^{\frac{1}{2}} q_{kl}}{\sigma_{kl}} - 2 w_{kl}q_{kl} \;\; \text{ with } \;\; \sigma_{kl} = \begin{cases}
a_{kl} &\text{if} \;\; l=1,\ldots,d_k,\\
b_k &\text{if} \;\; l=d_k+1,\ldots,B.
\end{cases}
\end{align*}

\begin{align*}
\text{Thus,}\;\;
    \nabla_{q_{kl}} \mathcal{L} = 0 & \iff W^{\frac{1}{2}} H_k W^{\frac{1}{2}} q_{kl} = \frac{w_{kl}\sigma_{kl}}{n_k} q_{kl}\,.
\end{align*}

It results that $q_{kl}$ is an eigenvector of $W^{\frac{1}{2}} H_k W^{\frac{1}{2}}$. The associated eigenvalue is $\lambda_{kl} = \frac{w_{kl}\sigma_{kl}}{n_k}$, which satisfies $\lambda_{kl} = q_{kl}^\prime W^{\frac{1}{2}} H_k W^{\frac{1}{2}} q_{kl}$.
Additionally, $q_{kl}^\prime q_{ku} = 0$ if $l \neq u$.

We also deduce from above that $\tr(W^{\frac{1}{2}} H_k W^{\frac{1}{2}}) = \sum_{l=1}^B \lambda_{kl}$, which implies that \begin{equation}\label{eqn:lambdas}\sum_{l=d_k+1}^B \lambda_{kl} = \tr(W^{\frac{1}{2}} H_k W^{\frac{1}{2}})  - \sum_{l=1}^{d_k} \lambda_{kl}.\end{equation}

Maximizing $\ell_{c_3}$ with respect to $q_{kl}$ is equivalent to minimizing $ -2 \ell_{c_3}$ with respect to $q_{kl}$.
Recall that
\begin{equation*} -2\ell_{c_3}(\textbf{c},\textbf{z},\textbf{v},\theta ) = \sum_{k=1}^K \left\{ n_k \left( \sum_{l=1}^{d_k} \log(a_{kl}) + \sum_{l=d_k +1}^B \log(b_{k})  \right) + n_k \left[ \sum_{l=1}^{d_k} \frac{\lambda_{kl}}{a_{kl}} + \sum_{l=d_k +1}^{B} \frac{\lambda_{kl}}{b_{k}}   \right]  \right\} + S
\end{equation*}
where $S$ does not depend of $(\lambda_{kl})$. 
Reporting expression \eqref{eqn:lambdas} in equation above gives:
\begin{align*}
-2\ell_{c_3}(\textbf{c},\textbf{z},\textbf{v},\theta )    &= \sum_{k=1}^K \left\{ n_k \left( \sum_{l=1}^{d_k} \log(a_{kl}) + \sum_{l=d_k +1}^B \log(b_{k})  \right)  \right.\\&   \left. \qquad \qquad + n_k  \sum_{l=1}^{d_k} \lambda_{kl} \left( \frac{1}{a_{kl}}  -  \frac{1 }{b_k} \right)   +  \frac{n_k}{b_k} tr(W^{\frac{1}{2}} H_k W^{\frac{1}{2}}) \right\} + S\,.
\end{align*}
Minimizing $ -2 \ell_{c_3}$ with respect to $q_{kl}$ is therefore equivalent to minimizing the quantity \[ \overline{\ell_{c_3}}(\textbf{c},\textbf{z},\textbf{v},\theta ) :=
\sum_{l=1}^{d_k} \lambda_{kl} \left( \frac{1}{a_{kl}}  - \frac{1}{b_k} \right).\] 

Remark that $\frac{1}{a_{kl}}  - \frac{1}{b_k} \leq 0$ (because $\forall \; l, a_{kl} \geq b_k$). As a consequence, $\overline{\ell_{c_3}}$ is an decreasing function of $\lambda_{kl}$. So the $l$th column $q_{kl}$ of the orthogonal matrix $Q_k$ is estimated by the eigenfunction associated to the $l^\text{th}$ highest eigenvalue of $W^{\frac{1}{2}} H_k W^{\frac{1}{2}}$.

\paragraph{Update of $a_{kl}$}
We first differentiate $\ell_{c_3}$ with respect to $a_{kl}$,
\begin{equation*}
    \frac{\partial\ell_{c_3}(\textbf{c},\textbf{z},\textbf{v},\theta )}{\partial a_{kl}} = -\frac{n_k}{2 a_{kl}} + \frac{n_k}{2} \frac{q_{kl}^\prime W^{\frac{1}{2}} H_k W^{\frac{1}{2}} q_{kl}}{a_{kl}^2} \end{equation*}
    Consequently,
    \[
    \frac{\partial\ell_{c_3}(\textbf{c},\textbf{z},\textbf{v},\theta )}{\partial a_{kl}} = 0 \iff a_{kl} =  q_{kl}^\prime W^{\frac{1}{2}} H_k W^{\frac{1}{2}} q_{kl}\,. 
\]

\begin{equation}
\text{Finally,}\;\;
    a_{kl}^{(m+1)} =  \lambda_{kl}^{(m)}, \;\; \text{where} \;\; \lambda_{kl}^{(m)} = q_{kl}^\prime W^{\frac{1}{2}} H_k^{(m)} W^{\frac{1}{2}} q_{kl}\,,\;\;  \text{and} \;\; \gamma_k^{(m)} = \sum_{i=1}^n t_{ik}^{(m)}
\end{equation}

\paragraph{Update of $b_k$}
Similarly,
\[
    \frac{\partial\ell_{c_3}(\textbf{c},\textbf{z},\textbf{v},\theta )}{\partial b_k} = - \frac{n_k}{2} \sum_{l=d_k+1}^B \frac{1}{b_k} + \frac{n_k}{2} \sum_{l=d_k+1}^B \frac{q_{kl}^\prime W^{\frac{1}{2}} H_k W^{\frac{1}{2}} q_{kl}}{b_k^2}. \]
    
    \begin{align*}
    \text{and thus,}\;\;
    \frac{\partial\ell_{c_3}(\textbf{c},\textbf{z},\textbf{v},\theta )}{\partial b_k} = 0 &\iff n_k \frac{B-d_k}{b_k} = \frac{n_k}{b_k^2} \sum_{l=d_k+1}^B q_{kl}^\prime W^{\frac{1}{2}} H_k W^{\frac{1}{2}} \\ 
    &\iff b_k = \frac{1}{B-d_k}\left( tr(W^{\frac{1}{2}} H_k W^{\frac{1}{2}}) - \sum_{l=1}^{d_k} \lambda_{kl}   \right)\,.
\end{align*}
\begin{equation}
    \text{We obtain that}\;\; b_k^{(m+1)} =  \frac{1}{B-d_k}\left( tr(W^{\frac{1}{2}} H_k^{(m)} W^{\frac{1}{2}}) - \sum_{l=1}^{d_k} \lambda_{kl}^{(m)}   \right)\,.
\end{equation}

\subsection{Second maximization step (CM2)}
At this step, we fix $\theta_1$ at $\theta_1^{(m+1)}$ and we compute $\theta^{(m+1)}$ as the value of $\theta_2$  that maximizes $Q(\theta\mid\theta^{(m)})$. 

Let recall the expressions of $\ell_{c_3}$: 
\begin{align*}
    \ell_{c_3}(\textbf{c},\textbf{z},\textbf{v},\theta )  =& \sum_{k=1}^K \sum_{i=1}^n z_{ik}v_{ik}\log g(c_i,\mu_k,\Sigma_k)\,-\frac{1}{2}\sum_{k=1}^K \sum_{i=1}^n z_{ik}(1-v_{ik})\left[B\log(2\pi)+\log|\Sigma_k|\right] - \\ &\,\frac{1}{2}\sum_{k=1}^K \sum_{i=1}^n z_{ik}(1-v_{ik})\left[    B \log (\eta_k) +   \eta_k^{-1} (c_i - \mu_k)^\prime\Sigma_k^{-1}(c_i-\mu_k)   \right]\,.
\end{align*}
By differentiating the expression $\ell_{c_3}$ with respect to $\eta_k$, we have:
\[
    \frac{\partial\ell_{c_3}(\textbf{c},\textbf{z},\textbf{v},\theta )}{\partial\eta_k} = \frac{-B}{\eta_k}\sum_{i=1}^n z_{ik}(1-v_{ik}) + \frac{1}{\eta_k^2}\sum_{i=1}^n z_{ik}(1-v_{ik}) (c_i - \mu_k)^\prime\Sigma_k^{-1}(c_i - \mu_k). \]
    We deduce that
 \begin{align*}   
   \frac{\partial\ell_{c_3}}{\partial\eta_k} = 0 &\iff \eta_k = \frac{\sum_{i=1}^n z_{ik}(1-v_{ik}) (c_i - \mu_k)^\prime\Sigma_k^{-1}(c_i - \mu_k)}{B\sum_{i=1}^n z_{ik}(1-v_{ik})}\,.
\end{align*}

Since $\eta_k > 1$, we rather set\begin{equation}
    \eta_k^{(m+1)} = \max\left\{ 1, \frac{\sum_{i=1}^n t_{ik}^{(m)}(1-q_{ik}^{(m)}) (c_i - \mu_k^{(m+1)}) ^\prime\Sigma_k^{-1(m+1)}(c_i - \mu_k^{(m+1)})}{B(\gamma_k^{(m)} - \gamma_{k_g}^{(m)})}  \right\}.
\end{equation}

\bibliographystyle{elsarticle-num} 
\bibliography{cas-refs}

\begin{thebibliography}{10}
\expandafter\ifx\csname url\endcsname\relax
  \def\url#1{\texttt{#1}}\fi
\expandafter\ifx\csname urlprefix\endcsname\relax\def\urlprefix{URL }\fi
\expandafter\ifx\csname href\endcsname\relax
  \def\href#1#2{#2} \def\path#1{#1}\fi

\bibitem{ferraty2003curves}
F.~Ferraty, P.~Vieu, Curves discrimination: a nonparametric functional
  approach, Computational Statistics \& Data Analysis 44~(1-2) (2003) 161--173.

\bibitem{ramsay2005springer}
J.~Ramsay, B.~Silverman, Functional data analysis, Springer Series in
  Statistics, 2005.

\bibitem{agyemang2006comprehensive}
M.~Agyemang, K.~Barker, R.~Alhajj, A comprehensive survey of numeric and
  symbolic outlier mining techniques, Intelligent Data Analysis 10~(6) (2006)
  521--538.

\bibitem{chandola2009anomaly}
V.~Chandola, A.~Banerjee, V.~Kumar, Anomaly detection: a survey, ACM Computing
  Surveys (CSUR) 41~(3) (2009) 1--58.

\bibitem{hodge2004survey}
V.~Hodge, J.~Austin, A survey of outlier detection methodologies, Artificial
  Intelligence Review 22~(2) (2004) 85--126.

\bibitem{chalapathy2019deep}
R.~Chalapathy, S.~Chawla, Deep learning for anomaly detection: a survey, arXiv
  preprint arXiv:1901.03407 (2019).

\bibitem{braei2020anomaly}
M.~Braei, S.~Wagner, Anomaly detection in univariate time-series: a survey on
  the state-of-the-art, arXiv preprint arXiv:2004.00433 (2020).

\bibitem{fraiman2001trimmed}
R.~Fraiman, G.~Muniz, Trimmed means for functional data, Test 10~(2) (2001)
  419--440.

\bibitem{cuevas2007robust}
A.~Cuevas, M.~Febrero, R.~Fraiman, Robust estimation and classification for
  functional data via projection-based depth notions, Computational Statistics
  22~(3) (2007) 481--496.

\bibitem{febrero2008outlier}
M.~Febrero, P.~Galeano, W.~Gonz{\'a}lez-Manteiga, Outlier detection in
  functional data by depth measures, with application to identify abnormal
  {NOx} levels, Environmetrics: The Official Journal of the International
  Environmetrics Society 19~(4) (2008) 331--345.

\bibitem{sun2011functional}
Y.~Sun, M.~G. Genton, Functional boxplots, Journal of Computational and
  Graphical Statistics 20~(2) (2011) 316--334.

\bibitem{hubert2015multivariate}
M.~Hubert, P.~J. Rousseeuw, P.~Segaert, Multivariate functional outlier
  detection, Statistical Methods \& Applications 24~(2) (2015) 177--202.

\bibitem{hubert2017multivariate}
M.~Hubert, P.~Rousseeuw, P.~Segaert, Multivariate and functional classification
  using depth and distance, Advances in Data Analysis and Classification 11~(3)
  (2017) 445--466.

\bibitem{dai2019directional}
W.~Dai, M.~G. Genton, Directional outlyingness for multivariate functional
  data, Computational Statistics \& Data Analysis 131 (2019) 50--65.

\bibitem{sakoe1978dynamic}
H.~Sakoe, S.~Chiba, Dynamic programming algorithm optimization for spoken word
  recognition, IEEE Transactions on Acoustics, Speech, and Signal Processing
  26~(1) (1978) 43--49.

\bibitem{giorgino2009computing}
T.~Giorgino, Computing and visualizing dynamic time warping alignments in {R}:
  the dtw package, Journal of Statistical Software 31~(7) (2009) 1--24.

\bibitem{sarda2017comparing}
A.~Sard{\'a}-Espinosa, Comparing time-series clustering algorithms in {R} using
  the dtwclust package, R package vignette 12 (2017) 41.

\bibitem{garcia2005proposal}
L.~A. Garcia-Escudero, A.~Gordaliza, A proposal for robust curve clustering,
  Journal of Classification 22~(2) (2005) 185--201.

\bibitem{cuesta1997trimmed}
J.~A. Cuesta-Albertos, A.~Gordaliza, C.~Matr{\'a}n, Trimmed $ k $-means: an
  attempt to robustify quantizers, The Annals of Statistics 25~(2) (1997)
  553--576.

\bibitem{staerman2019functional}
G.~Staerman, P.~Mozharovskyi, S.~Cl{\'e}men{\c{c}}on, F.~d’Alch{\'e} Buc,
  Functional isolation forest, in: Asian Conference on Machine Learning, PMLR,
  2019, pp. 332--347.

\bibitem{liu2008isolation}
F.~T. Liu, K.~M. Ting, Z.-H. Zhou, Isolation forest, in: 2008 Eighth IEEE
  International Conference on Data Mining, IEEE, 2008, pp. 413--422.

\bibitem{jacques2014functional}
J.~Jacques, C.~Preda, Functional data clustering: a survey, Advances in Data
  Analysis and Classification 8~(3) (2014) 231--255.

\bibitem{abraham2003unsupervised}
C.~Abraham, P.-A. Cornillon, E.~Matzner-L{\o}ber, N.~Molinari, Unsupervised
  curve clustering using {B}-splines, Scandinavian Journal of Statistics 30~(3)
  (2003) 581--595.

\bibitem{peng2008distance}
J.~Peng, H.-G. M{\"u}ller, Distance-based clustering of sparsely observed
  stochastic processes, with applications to online auctions, The Annals of
  Applied Statistics 2~(3) (2008) 1056--1077.

\bibitem{ieva2011multivariate}
F.~Ieva, A.~Paganoni, D.~Pigoli, V.~Vitelli, Multivariate functional clustering
  for the analysis of {ECG} curves morphology, in: Cladag 2011 (8th
  International Meeting of the Classification and Data Analysis Group), 2011,
  pp. 1--4.

\bibitem{james2003clustering}
G.~M. James, C.~A. Sugar, Clustering for sparsely sampled functional data,
  Journal of the American Statistical Association 98~(462) (2003) 397--408.

\bibitem{heard2006quantitative}
N.~A. Heard, C.~C. Holmes, D.~A. Stephens, A quantitative study of gene
  regulation involved in the immune response of anopheline mosquitoes: an
  application of bayesian hierarchical clustering of curves, Journal of the
  American Statistical Association 101~(473) (2006) 18--29.

\bibitem{bouveyron2011model}
C.~Bouveyron, J.~Jacques, Model-based clustering of time series in
  group-specific functional subspaces, Advances in Data Analysis and
  Classification 5~(4) (2011) 281--300.

\bibitem{jacques2014model}
J.~Jacques, C.~Preda, Model-based clustering for multivariate functional data,
  Computational Statistics \& Data Analysis 71 (2014) 92--106.

\bibitem{jacques2013funclust}
J.~Jacques, C.~Preda, Funclust: a curves clustering method using functional
  random variables density approximation, Neurocomputing 112 (2013) 164--171.

\bibitem{schmutz2020clustering}
A.~Schmutz, J.~Jacques, C.~Bouveyron, L.~Cheze, P.~Martin, Clustering
  multivariate functional data in group-specific functional subspaces,
  Computational Statistics (2020) 1--31.

\bibitem{punzo2016parsimonious}
A.~Punzo, P.~D. McNicholas, Parsimonious mixtures of multivariate contaminated
  normal distributions, Biometrical Journal 58~(6) (2016) 1506--1537.

\bibitem{tomarchio2020dichotomous}
S.~D. Tomarchio, A.~Punzo, Dichotomous unimodal compound models: application to
  the distribution of insurance losses, Journal of Applied Statistics
  47~(13-15) (2020) 2328--2353.

\bibitem{punzo2019new}
A.~Punzo, A new look at the inverse {G}aussian distribution with applications
  to insurance and economic data, Journal of Applied Statistics 46~(7) (2019)
  1260--1287.

\bibitem{punzo2018fitting}
A.~Punzo, A.~Mazza, A.~Maruotti, Fitting insurance and economic data with
  outliers: a flexible approach based on finite mixtures of contaminated gamma
  distributions, Journal of Applied Statistics 45~(14) (2018) 2563--2584.

\bibitem{delaigle2010defining}
A.~Delaigle, P.~Hall, Defining probability density for a distribution of random
  functions, The Annals of Statistics 38~(2) (2010) 1171--1193.

\bibitem{yakowitz1968identifiability}
S.~J. Yakowitz, J.~D. Spragins, On the identifiability of finite mixtures, The
  Annals of Mathematical Statistics (1968) 209--214.

\bibitem{di2007mixture}
M.~Di~Zio, U.~Guarnera, R.~Rocci, A mixture of mixture models for a
  classification problem: the unity measure error, Computational Statistics \&
  Data Analysis 51~(5) (2007) 2573--2585.

\bibitem{punzo2016contaminatedmixt}
A.~Punzo, A.~Mazza, P.~D. McNicholas, Contaminated{M}ixt: an {R} package for
  fitting parsimonious mixtures of multivariate contaminated normal
  distributions, Journal of Statistical Software 85 (2018).

\bibitem{dempster1977maximum}
A.~P. Dempster, N.~M. Laird, D.~B. Rubin, Maximum likelihood from incomplete
  data via the {EM} algorithm, Journal of the Royal Statistical Society: Series
  B (Methodological) 39~(1) (1977) 1--22.

\bibitem{meng1993maximum}
X.-L. Meng, D.~B. Rubin, Maximum likelihood estimation via the {ECM} algorithm:
  a general framework, Biometrika 80~(2) (1993) 267--278.

\bibitem{browne2011model}
R.~P. Browne, P.~D. McNicholas, M.~D. Sparling, Model-based learning using a
  mixture of mixtures of {G}aussian and uniform distributions, IEEE
  Transactions on Pattern Analysis and Machine Intelligence 34~(4) (2011)
  814--817.

\bibitem{hartigan1975clustering}
J.~A. Hartigan, Clustering algorithms, John Wiley \& Sons, Inc., 1975.

\bibitem{cattell1966scree}
R.~B. Cattell, The scree test for the number of factors, Multivariate
  Behavioral Research 1~(2) (1966) 245--276.

\bibitem{biernacki2000assessing}
C.~Biernacki, G.~Celeux, G.~Govaert, Assessing a mixture model for clustering
  with the integrated completed likelihood, IEEE Transactions on Pattern
  Analysis and Machine Intelligence 22~(7) (2000) 719--725.

\bibitem{birge2007minimal}
L.~Birg{\'e}, P.~Massart, Minimal penalties for {G}aussian model selection,
  Probability Theory and Related Fields 138~(1-2) (2007) 33--73.

\bibitem{akaike1974new}
H.~Akaike, A new look at the statistical model identification, IEEE
  Transactions on Automatic Control 19~(6) (1974) 716--723.

\bibitem{schwarz1978estimating}
G.~Schwarz, Estimating the dimension of a model, The Annals of Statistics 6~(2)
  (1978) 461--464.

\bibitem{fraley1998many}
C.~Fraley, A.~E. Raftery, How many clusters? {W}hich clustering method?
  {A}nswers via model-based cluster analysis, The Computer Journal 41~(8)
  (1998) 578--588.

\bibitem{bouveyron2015discriminative}
C.~Bouveyron, E.~C{\^o}me, J.~Jacques, The discriminative functional mixture
  model for a comparative analysis of bike sharing systems, The Annals of
  Applied Statistics 9~(4) (2015) 1726--1760.

\bibitem{preda2007regression}
C.~Preda, Regression models for functional data by reproducing kernel {H}ilbert
  spaces methods, Journal of Statistical Planning and Inference 137~(3) (2007)
  829--840.

\bibitem{rand1971objective}
W.~M. Rand, Objective criteria for the evaluation of clustering methods,
  Journal of the American Statistical Association 66~(336) (1971) 846--850.

\bibitem{forbes2014new}
F.~Forbes, D.~Wraith, A new family of multivariate heavy-tailed distributions
  with variable marginal amounts of tailweight: application to robust
  clustering, Statistics and Computing 24~(6) (2014) 971--984.

\bibitem{morris2019asymmetric}
K.~Morris, A.~Punzo, P.~D. McNicholas, R.~P. Browne, Asymmetric clusters and
  outliers: mixtures of multivariate contaminated shifted asymmetric {L}aplace
  distributions, Computational Statistics \& Data Analysis 132 (2019) 145--166.

\bibitem{punzo2021multiple}
A.~Punzo, C.~Tortora, Multiple scaled contaminated normal distribution and its
  application in clustering, Statistical Modelling 21~(4) (2021) 332--358.

\bibitem{bellas2013model}
A.~Bellas, C.~Bouveyron, M.~Cottrell, J.~Lacaille, Model-based clustering of
  high-dimensional data streams with online mixture of probabilistic {PCA},
  Advances in Data Analysis and Classification 7~(3) (2013) 281--300.

\end{thebibliography}

\end{document}